\documentclass[11pt]{article}

\usepackage{amsmath,amsfonts,amsthm,amssymb,amscd,color}
\renewcommand{\Re}{\mathop{\rm Re}\nolimits} \renewcommand{\Im}{\mathop{\rm
Im}\nolimits}

\theoremstyle{plain}
\newtheorem{theorem}{Theorem}[section]
\newtheorem{lemma}[theorem]{Lemma}
\newtheorem{proposition}[theorem]{Proposition}
\newtheorem{corollary}[theorem]{Corollary} \theoremstyle{definition}
 \theoremstyle{remark}
\newtheorem{remark}[theorem]{Remark} 
 
\newcommand{\R}{{\mathbb R}}

 \def\im{{\rm i}}
\newcommand{\C}{\mathbb{C}}

\def\({\left(}
\def\){\right)}
\def\<{\left\langle}
\def\>{\right\rangle}

\numberwithin{equation}{section}

\begin{document}

\title{\bf  The asymptotic stability of solitons in the cubic NLS equation on the line}

\author{Scipio Cuccagna\footnote{Department of Mathematics and Geosciences,  University
of Trieste, Trieste, 34127 Italy, e-mail: scuccagna@units.it} \hspace{0.1cm} and
Dmitry E. Pelinovsky\footnote{Department of Mathematics and Statistics, McMaster
University, Hamilton ON, Canada, L8S 4K1, e-mail: dmpeli@math.mcmaster.ca}}

\maketitle

\begin{abstract}
We use the inverse scattering transform, the auto--B\"{a}cklund transformation
and the steepest descent method  of Deift and Zhou to obtain the asymptotic stability of the solitons
in the cubic NLS (nonlinear Schr\"{o}dinger) equation.
\end{abstract}

\section{Introduction}

We consider the Cauchy problem for the cubic focusing NLS (nonlinear Schr\"odinger)
equation on $\R$:
\begin{equation}\label{eq:nls}
 \im u_{t }+  u_{xx } + 2  |u|^2  u=0, \quad u(0)=u_0.
\end{equation}
The Cauchy problem \eqref{eq:nls} is globally well posed in $L^{2}(\R  )$, by
the following result due to Tsutsumi \cite{tsutsumi}.

\begin{theorem}
\label{th:L2} Given $u_0\in L^{2}(\R)$, then there exists a unique
solution
$$
u(t) \in C^0 (\R, L^2(\R ))\cap L^4_{loc}(\R, L^\infty(\R ))
$$
of the integral equation
\begin{equation}\label{eq:nls1}
u(t) = e ^{\im t\partial _x^2 }u_0 +2\im \int _0^t  e ^{\im (t-t' )\partial _x^2 }  |u(t' )|^2u(t' ) dt'  .
\end{equation}
We have $\| u(t)\| _{L^2}=\| u_0\| _{L^2}$.
Furthermore, if $ u_{0n} \rightarrow u_{0 } $ as $n \to \infty$
in  $  L^{2 } (\R ) $ and $u_n(t)$ is the solution of the NLS equation
with $u_n(0) = u_{0n}$, then,   as $n \to \infty$, for any $t \in \R$ we have
$ u_{ n}(t) \rightarrow u (t) $  in  $  L^{2 } (\R ) $.
\end{theorem}

\begin{remark}
An important class of solutions of the NLS equation are the solitons, defined by
\begin{equation}\label{eq:groundstate}
\varphi_{\omega,\gamma,v}(t,x-x_0) :=  \omega  e^{\im xv +
\im ( \omega ^2 - v^2)t + \im \gamma}\  \text{sech} ( \omega (x-2 vt - x_0)).
\end{equation}
We are interested here to the question of their asymptotic stability, when
$u_0$ is close to $\varphi_{\omega,\gamma,v}$ for a particular $(\omega,\gamma,v)$.
\end{remark}

\noindent {\em Notations:} the following Hilbert spaces are  the closures
of  the space $C_0^\infty(\mathbb{ R})$ with respect to the following norms,
where  $ \langle x \rangle := \sqrt{ 1+|x|^2} $:
\begin{equation*}
\begin{aligned} & L^{2 ,s}(\mathbb{ R}  ) \text{  defined with  }
  \| u \| _{L^{2,s}(\R) }  :=  \|   \langle x \rangle ^s   u \| _{L^2 (\mathbb{ R}  )} ;  \\
  &   \dot H^s(\mathbb{ R}  ) \text{  defined with  }
 \| u \| _{\dot H^s (\mathbb{ R}  ) }  :=  \|   |x | ^s   \widehat{u} \| _{L^2 (\mathbb{ R}  )}
  \text{   where   $ \widehat{u}$ is the Fourier transform;} \\
 &  H^s(\mathbb{ R}  ) \text{  defined with  }
  \| u \| _{H^s(\mathbb{ R}  ) }  :=  \|   \langle x \rangle ^s   \widehat{u} \| _{L^2 (\mathbb{ R}  )},
  \text{ such that } H^s(\mathbb{ R}  )=\dot H^s (\mathbb{ R}  )\cap L^2 (\mathbb{ R}).
 \end{aligned}
\end{equation*}
We set $\Sigma _s:=H^s(\mathbb{ R}  )\cap  L^{2 ,s}(\mathbb{ R}  )$.
  Our aim  is to prove the following result.

\begin{theorem}\label{th:asstab}  Fix $s\in (1/2, 1]$. Consider the NLS soliton
$\varphi_{\omega_0,\gamma_0,v_0}(0,x-x_0)$.
Then, there exist positive constants
$\varepsilon _0=\varepsilon _0(\omega _0, v_0)$,  $T=T(\omega _0, v_0)$  and $C=C(\omega _0, v_0)$
such that if $u_0\in  L^{2,s}(\R )$ and if
\begin{equation}\label{eq:asstab0}
\begin{aligned} & \epsilon :=\| \varphi  _{\omega _0, \gamma _0, v_0}(0,\cdot -x_0) -
u_0\|  _{L^{2,s}(\R)} <  \varepsilon _0 ,
 \end{aligned}
\end{equation}
then there exist  two   ground states
$\varphi_{\omega_1 ,\gamma_\pm ,v_1}(t,x-x_\pm )$
such that for the solution of the Cauchy problem \eqref{eq:nls}   provided by Theorem \ref{th:L2}, we have
\begin{equation}\label{eq:asstab1}
| (\omega _1, \gamma _\pm , v_1 , x_\pm )- (\omega _0, \gamma _0, v_0 , x_0)|< C\epsilon
\end{equation}
and, for all $\pm t \geq T$,
\begin{equation}
\label{eq:asstab2}
\|  u(t,\cdot ) - \varphi _{\omega _1, \gamma _\pm , v_1}(t,\cdot -x_\pm )\| _{L^\infty  (\R )}<
C \epsilon |t|^{-\frac{1}{2}}.
\end{equation}
\end{theorem}

In general the two   ground states $\varphi_{\omega_1 ,\gamma_\pm ,v_1}(t,x-x_\pm )$
are distinct, see  {Lemma}
  \ref{lem:distgs} at the end of the paper.

A key ingredient  for the above result comes from the methods of the inverse scattering transform (IST)
theory, found in references \cite{BC,Z,DZ2,DZ,DM}. In particular, we use
the steepest descent method and the auto--B\"{a}cklund transformation
discussed in  \cite{DP}. Theorem \ref{th:asstab} is an analogue of the results
about the asymptotic behavior of solutions decaying to 0, obtained in \cite{DZ2,HN,DZ,DM,DP}.
Compared to these references, we do not reproduce in Theorem \ref{th:asstab}
the asymptotic expansions of the solution $u$ for large values of $t$, but we ease
the restrictions on the initial data by allowing $u_0\in  L^{2,s}(\R )$ for  $s\in (1/2, 1]$
and not just for $s=1.$

Theorem \ref{th:asstab} should be contrasted to the results for non-integrable systems,
where the orbits of the  solitons   which attract the solution $u(t)$
are  presumably  not the same as $t\to +\infty$ and $t\to -\infty$, see   \cite{SW2,BP1,PiW}
for early results. In  the case of the cubic NLS equation, it turns out that the selected
asymptotic soliton is simply defined by the eigenvalue of a spectral problem
supported by the initial datum $u_0$  but it has a different reflection  coefficient, which
is zero for the solitons (\ref{eq:groundstate}) and nonzero for a generic $u_0$.

Another feature of  non-integrable systems, is that  the rate of decay
in the right hand side of \eqref{eq:asstab2} is generally slower, because of metastable states
which are not present for the cubic NLS equation. The theory how to treat these
metastable states was initiated in  \cite{BP2,SW3} and for recent developments and
further references  we refer to \cite{Cu2,Cu3,Cu0}.   Obviously the absence of
metastable states for the cubic NLS equation simplifies the discussion.
Notice that  \cite{Kevrekidis} conjectures   the non-existence of metastable states
in integrable systems.

Theorem  \ref{th:asstab} appears to be out of reach of the perturbative methods initiated in
\cite{SW2,BP1,PiW} and developed in  a number of papers using a similar framework.
This is because of the "strength"  of the cubic nonlinearity in the cubic NLS equation.
This strength is   responsible for the fact that the classical result in
\cite{MS} on the dispersion of small solutions   of $L^2$ subcritical equations
does not apply to the cubic NLS equation, although  it was proved also for the  cubic NLS equation    a decade later in \cite{HN},
with an approach    similar  to \cite{MS} but with an additional normal form argument.
The results in  \cite{MS,HN} are based on {\it invariant fields}  which exploit
symmetries of the equations not present in the case of the linearization of an NLS at a soliton. And while
\cite{MS} has been partially extended to settings without translation symmetry in \cite{CGV},
so far   the approach in \cite{MS,HN} has never been applied
directly to the problem of asymptotic stability of the solitons.

Mizumachi and Pelinovsky \cite{MP} proposed to treat  the orbital stability of solitons of  the cubic NLS
equation by using an auto--B\"{a}cklund transformation which transforms a soliton in the zero
solution and preserves the equation. They  proved that the B\"{a}cklund  transformation
is a homeomorphism  in $L^2$.   The  B\"{a}cklund  transformation can then be used to
transfer  Theorem \ref{th:L2} into a statement about solutions close to the soliton in  $L^2(\R )$,
in particular proving that
solitons of the cubic NLS equation are  orbitally stable in $L^2(\R )$,
thus transposing to
$L^2(\R )$ the classical result of orbital  stability proved for the space $H^1(\R )$ in \cite{W1}.
In \cite{MP}, a discussion was initiated on the possible use of the same transformation
to transfer the  dispersion scattering result for small  solutions  in \cite{HN} to an asymptotic stability
result for solitons in $\Sigma_1$.   However  it is an open   question   whether or not the  B\"{a}cklund  transformation
in \cite{MP}   is a homeomorphism   in $\Sigma _1$.

The inspiration for the present paper comes  however from a   paper by Deift and Park \cite{DP},
where there is a particularly simple and   explicit      B\"{a}cklund  transformation,
see  \eqref{eq:darboux2} later.   Using the  steepest descent method of Deift and Zhou \cite{DZ3}
it is possible to bound all the terms
of formula    \eqref{eq:darboux2} and prove  Theorem \ref{th:asstab}. Specifically,
  by means of direct scattering,
it is possible to derive the spectral data associated with the solution  $u$ of the
Cauchy problem (\ref{eq:nls}).  Then, from   mapping properties of the
inverse scattering transform proved in \cite{Z,DZ}, which are similar to   mapping properties of the
inverse Fourier transform, the solution  $u$ is expressed by means of the transformation
formula \eqref{eq:darboux2} as the sum of a pure radiation solution  $\widetilde{u}$
and an appropriate fraction of Jost functions associated to the potential  $\widetilde{u}$.
The results in \cite{DZ,DM,DP} are applied directly to the pure radiation solution $\widetilde{u}$.
Also Jost functions and their fraction can be easily analyzed using other results from
\cite{DZ,DM,DP}. This yields Theorem \ref{th:asstab}. Notice that Theorem \ref{th:purerad} in Sect.\ref{sec:purerad}   extends the result in \cite{DZ,DM} to initial data
in $L^{2,s}(\R )$ for    all $s\in (1/2, 1]$.

We do not make any particular claim of originality, since Theorem \ref{th:asstab}
is a natural corollary of the previous works \cite{DZ,DM,GT,DP}. Nonetheless, we feel that
it is important  that  Theorem \ref{th:asstab} be stated  explicitly and proved.

The paper is organized as follows. Section \ref{sec:scatt} gives details of the
direct and inverse scattering transforms for the cubic NLS equation. Section \ref{sec:purerad}
contains a review of the asymptotic scattering theory for the pure radiation solution.
Section \ref{sec:asstab} explains the arguments needed to prove the asymptotic stability of solitons
formulated in Theorem \ref{th:asstab}.

{\bf Acknowledgements.} S.C. was partially funded  by a  grant  FRA 2009
from the University of Trieste and by the grant FIRB 2012 (Dinamiche Dispersive).
D.P. was partially funded by the NSERC Discovery grant.

\section{Direct and Inverse Scattering Transforms}
\label{sec:scatt}

The Cauchy problem (\ref{eq:nls}) for the cubic NLS equation can be solved
through  the direct and inverse scattering transform.

Consider a function  $u(x) \in L^1(\R)$ and recall
that $L^{2,s}(\R )$ is embedded into $L^{1 } (\R )$ for any $s > \frac{1}{2}$.
The spectral system associated with the cubic NLS equation takes the form:
\begin{equation}  \label{eq:zs}
\psi _x     = -\im    z \sigma _3 \psi + Q(u(x)) \psi,
\end{equation}
where
\begin{equation*}
Q(u(x)) := \begin{pmatrix} 0  &  u( x) \\ - \overline{u}( x)
  & 0 \end{pmatrix} \, , \quad   \sigma _3 :=    \begin{pmatrix}
1  & 0 \\
0 & -1 \end{pmatrix}.
\end{equation*}

Set $e_1= (1,0)^T $  and  $e_2= (0,1)^T $. According to the
direct scattering theory \cite{ABT}, for any fixed $ z \in   \C _+$ (i.e. $\Im z >0$) there exists a unique
$ \C ^2$ valued  solution $\phi   ( x,z) $ of
the spectral system \eqref{eq:zs} such that
\begin{equation}
\label{function-phi}
\lim _{x\to -\infty } \phi   ( x,z) e^{  \im x z} = e_1  \quad \text{ and }\quad
\lim _{x\to +\infty } \phi   ( x,z) e^{  \im x z} = a( z)e_1,
\end{equation}
where  $a(z)$ is an analytic function in $\C _+$, continuous in   $\overline{\C} _+$ with $\lim _{z\to \infty}a( z)=1$.
We call $a(z)$  the {\it scattering function}. The following
result is well known (see, i.e., \cite{ABT,BC}).

\begin{lemma}
\label{lem:dirsc1}
There exists  an open dense set $\mathcal{G} \subset L^{ 1}(\R )$   such that, for $u \in \mathcal{G}$,
the scattering function $a(z) $ has at most a  finite number of
zeros  forming a set $\mathcal{Z}_+=\{ z_1,..., z_n\} $ in $ {\C} _+$,
with   $a(z)\neq 0$ for all $z\in \R$ and $a'(z_k)\neq 0 $ for all $k$.
The cardinality $u\to \sharp \mathcal{Z}_+ $ is locally constant
near $u$ in $\mathcal{G}$ and  the map
$\mathcal{G} \ni u \to (z_1,...,z_n) \in \C_+^n$ is locally Lipschitz.
\end{lemma}
We denote by $\mathcal{G}_n$ the open subset of $\mathcal{G}$ formed
by the elements such that $\sharp \mathcal{Z}_+ = n$.
\begin{remark}\label{rem-soliton-generic}
We call the potentials in  $\mathcal{G}$  \textit{generic}.   $\mathcal{G}_1$ contains
the solitons (\ref{eq:groundstate}). Notice  that
  small $L^1$-perturbations to the solitons are in $\mathcal{G}_1$.
     See \cite{PelSul}
for  other integrable equations  where this is not true.
\end{remark}

There exists a  unique $\C ^2$ valued solution $\psi( x,z) $ of
the spectral system \eqref{eq:zs} satisfying
\begin{equation}
\label{function-psi}
\lim _{x\to +\infty } \psi( x,z) e^{-\im x z} = e_2.
\end{equation}
If $u \in \mathcal {G}$ and $\mathcal{Z}_+\neq \emptyset$,
then for each $z_k\in  \mathcal{Z}_+$ we have $\phi( x,z_k)=\gamma _k  \psi   ( x,z_k) $
for some $\gamma _k \in \C _* :=\C \backslash \{ 0\}$.
Set $c_k =\gamma_k / a'(z_k )$ and call it {\it the norming constant}.

For $z\in \R$, the solution of  the spectral system \eqref{eq:zs}
with the boundary value
\begin{equation}
\lim _{x\to -\infty } \phi   ( x,z) e^{  \im x z} = e_1
\end{equation}
satisfies the scattering problem
\begin{equation}
\lim _{x\to +\infty } \left[ \phi   ( x,z) e^{  \im x z} -a( z) e_1 -e^{2  \im x z}b(   z) e_2
\right] = 0,
\label{scattering-problem}
\end{equation}
where $b(z)$ is a continuous function on $\R$. Set
$r(z):=\frac{b(   z)}{a( z)}$ and call it the {\it reflection coefficient}.

We consider now the Jost functions  defined by
the Volterra integral equations, see \cite{ABT},
\begin{equation}\label{eq:jost}
\begin{aligned}
  &   \mathfrak{m} _1^{\pm }(x,z) = e_1   +\int _{\pm \infty}^x   \begin{pmatrix}
    1   & 0 \\  0
  & e^{2\im  (x-y) z }
 \end{pmatrix}    Q(u(y))  \mathfrak{m} _1^{\pm } (y,z)  dy  ,\\&
 \mathfrak{m} _2^{\pm } (x,z) = e_2  +\int _{\pm \infty}^x     \begin{pmatrix}
    e^{-2\im  (x-y) z }   & 0 \\  0
  & 1
 \end{pmatrix}    Q(u(y))   \mathfrak{m} _2^{\pm }(y,z)  dy  .
\end{aligned}
\end{equation}
The functions $\mathfrak{m}_1^-(x,z)$ and $\mathfrak{m}_2^+(x,z)$ are analytic for $z \in \mathbb{C}_+$,
whereas the functions $\mathfrak{m}_2^-(x,z)$ and $\mathfrak{m}_1^+(x,z)$ are analytic for $z \in \mathbb{C}_-$,
\cite{ABT}.

\begin{remark}
In terms of functions $\phi$ and $\psi$ introduced in (\ref{function-phi}) and (\ref{function-psi}),
we have $\mathfrak{m}_1^-(x,z) = \phi(x,z) e^{i x z}$
and $\mathfrak{m}_2^+(x,z) = \psi(x,z) e^{-\im x z}$ for $z \in \C _+$.
\end{remark}

It follows from the scattering problem (\ref{scattering-problem})
and the Wronskian identities for the spectral system (\ref{eq:zs}),
we have for $z\in \R$,
\begin{equation}
a(z)=\det[\mathfrak{m} _1^{- }(x,z) , \mathfrak{m} _2^{+ }(x,z)]
\end{equation}
and
\begin{equation}
b(z)=\det[\mathfrak{m} _1^{+ }(x,z) , e^{-2\im xz}\mathfrak{m} _1^{- }(x,z)],
\end{equation}
where matrices $[\cdot,\cdot]$ are defined in the sense of column vectors
and the Wronskian determinants are $x$-independent. The following result is obtained
with a minor modification of the argument in Theorem 3.2   \cite{DZ}.

\begin{lemma}
  \label{lem:dirsc2}
  Let  $s\in \left(\frac{1}{2},1\right]$. For  $ u \in  L^{2,s} (\R )\cap \mathcal{G} $ we have
  $r\in H^s(\R) $. Furthermore, the map $L^{2,s} (\R )\cap \mathcal{G} \ni u \to r \in H^s(\R)$
  is locally Lipschitz.
\end{lemma}
\proof
We make the following claim:  for any fixed $\kappa _0>0$ there exists
a positive constant $C$ such that if  $\| u \| _{L^{2,s}(\R)} \le \kappa _0$,
then we have for $j = 1,2$:
\begin{equation} \label{eq:jost0}
\begin{aligned}
  &      \|    \mathfrak{m} _j^{- }(x,z )  -e_j \|   _{   H^s_z (\R )  }
\le C     \| u \| _{L^{2,s}(\R)}   \quad \text{  for all $x\le 0$}\\&   \|    \mathfrak{m} _j^{+ }(x,z )  -e_j \|   _{   H^s_z (\R )   }
\le C   \| u \| _{L^{2,s}(\R)}   \quad \text{  for all $x\ge 0$}.
\end{aligned}
\end{equation}
Let us assume \eqref{eq:jost0} for a moment. Then $b\in  H^s ( \R) $ because
\begin{eqnarray}
\nonumber
  b(z) & = & \det\left[\mathfrak{m} _1^{+ }(0,z) ,  \mathfrak{m} _1^{- }(0,z)\right] \\
\nonumber
  & = & \det\left[\mathfrak{m} _1^{+ }(0,z)-e_1 ,  \mathfrak{m} _1^{- }(0,z)-e_1\right] \\
  & \phantom{t} & \phantom{texttext} +
  \det\left[\mathfrak{m} _1^{+ }(0,z)-e_1 ,   e_1\right] +  \det\left[e_1 ,  \mathfrak{m} _1^{- }(0,z)-e_1\right],
    \label{b-bound}
\end{eqnarray}
where we recall that $H^s(\R)$ is a Banach algebra with respect to pointwise multiplication
for any $s > \frac{1}{2}$. Similarly,  $(a-1)\in  H^s ( \R) $ because
\begin{eqnarray}
\nonumber
a(z) & = & \det\left[\mathfrak{m} _1^{- }(0,z) ,  \mathfrak{m} _2^{+ }(0,z)\right] \\
\nonumber
& = &  \det\left[\mathfrak{m} _1^{- }(0,z)-e_1 ,  \mathfrak{m} _2^{+ }(0,z)\right]
+  \det\left[e_1 ,  \mathfrak{m} _2^{+ }(0,z)-e_2\right]   +  \det\left[e_1 ,   e_2\right] \\
\nonumber
& = & 1 + \det\left[\mathfrak{m} _1^{- }(0,z)-e_1 ,  \mathfrak{m} _2^{+ }(0,z)-e_2\right] \\
& \phantom{t} & \phantom{texttext}
  +  \det\left[e_1 ,  \mathfrak{m} _2^{+ }(0,z)-e_2\right] +  \det\left[\mathfrak{m} _1^{- }(0,z)-e_1,e_2\right].
  \label{a-bound}
\end{eqnarray}
We conclude that  if  $u\in  L^{2,s}(\R) \cap \mathcal{G}$ then $r\in H^s(\R)$.
So this shows that  we have a map
$L^{2,s}(\R) \cap \mathcal{G} \ni u \to r \in H^s(\R)$.
We skip the proof of the fact that the map
$L^{2,s}(\R) \cap \mathcal{G} \ni u \to r \in H^s(\R)$ is locally Lipschitz.

We now prove \eqref{eq:jost0}. It suffices to consider the case $j=1$ and the minus sign only.
The proof is   based on the fact that if  there is $s \in (0,1]$ such that for an $f\in L^2(\R )$ we have
  \begin{equation}\label{eq:sob}
\|  f(\cdot +h) -f(\cdot )\| _{L^2(\R )}\le C |h|^s, \quad \forall h\in \R,
\end{equation}
then $f\in H^s(\R )= \dot H^s(\R )\cap L^2(\R )$ and there is a  positive constant $c$  independent of $f$,
such that $\| f  \| _{ \dot H^s(\R) }\le c\  C$.

Let us define
\begin{equation*}
K f(x,z) := \int _{- \infty}^x    \begin{pmatrix}
    1   & 0 \\  0  & e^{2\im  (x-y) z } \end{pmatrix}   Q(u(y))  f(y,z)  dy.
\end{equation*}
By Theorem 3.2 \cite{DZ} (see also \cite{ABT}), we have
\begin{equation*}
\| K   e_2\| _{L^\infty _x (\R ,  L^2_z (\R )) }\le \| u \|_{L^2(\R )}
\end{equation*}
and for any $x_0 \le +\infty$,
\begin{equation*}
\left\|   (1-K )^{-1} \right\|_{L^\infty _x ((-\infty,x_0) ,  L^2_z (\R )) \to  L^\infty_x ((-\infty,x_0),L^2_z (\R ))}
\le e^{\sqrt{2}\| u \| _{L^1}}.
\end{equation*}
Furthermore, for $x\le 0$ we have
 \begin{equation}\label{eq:jost-2}\begin{aligned}
\|   \mathfrak{m} _1^{- }(x,z ) -e_2\|   _{   L^2_z (\R )  }   \le &
  e^{\sqrt{2}\| u \| _{L^1}}  \left\|  \int _{- \infty}^x e^{2\im  (x-y) z } \overline{u}(y) dy
  \right\|   _{   L^2_z (\R )  } \\
    \le &   e^{\sqrt{2}\| u \| _{L^1}}
  \left (  \int _{- \infty}^x  \langle y \rangle ^{2s} |{u}(y)|^2  dy    \right ) ^{1/2} \langle x \rangle ^{-s}  \\
    \le &
  e^{\sqrt{2}\| u \| _{L^1}} \| u \| _{L^{2,s}(\R)}  \langle x \rangle ^{-s}  .\end{aligned}
\end{equation}
To complete the proof of \eqref{eq:jost0} for $j=1$ and the  minus sign it is enough to prove
and estimate of the form \eqref{eq:sob}  with $C\lesssim  \| u \| _{L^{2,s}(\R)} $.
Define $n(x,z ):=   \mathfrak{m} _1^{- }(x,z+h)   -  \mathfrak{m} _1^{- }(x,z ) $ for $h\in \R$. We have
\begin{eqnarray}
\nonumber
(1-K )  n (x,z ) & = & \int _{- \infty}^x   \begin{pmatrix}
    0   & 0 \\  0  & e^{2\im  (x-y)  (z+h) }  -e^{2\im  (x-y)  z }
 \end{pmatrix}    Q(u(y))  (\mathfrak{m} _1^{- } (y,z)-e_1)  dy \\
 & \phantom{t} & \phantom{texttext} +  \int _{- \infty}^x   \begin{pmatrix}
    0      \\   \left ( e^{2\im  (x-y)  (z+h) }  -e^{2\im  (x-y)  z } \right ) \overline{u}(y)
 \end{pmatrix}     dy. \label{eq:jost1}
\end{eqnarray}
Using the Fourier transform ${\mathcal F}$, we have for $x\le 0$,
\begin{equation} \label{eq:jost2}
\begin{aligned}
  & \left\|   \int _{- \infty}^x \left  ( e^{2\im  (x-y)  (z+h) }  -e^{2\im  (x-y)  z } \right ) \overline{u}(y)
      dy  \right\|  _{L^2_z }  \\& =  \|    {\mathcal F}^*[u (\cdot +x) \chi _{\R _-} ] (z+h)-  {\mathcal F}^*[u (\cdot +x) \chi _{\R _-} ] (z )  \|  _{L^2_z } \\& \le C \|      {\mathcal F}^*[u (\cdot +x) \chi _{\R _-} ] (z )  \|  _{H^s_z(\R) } |h|^s = \|   u (y +x)   \|  _{L^{2,s}_y(\R _-) } |h|^s \le  \|   u     \|  _{L^{2,s}(\R)  } |h|^s
\end{aligned}
\end{equation}
and, using estimate \eqref{eq:jost-2},
\begin{align}  \label{eq:jost3}
  & \|    \text{first term r.h.s. \eqref{eq:jost1}} \|  _{L^2_z }  \le    2^{1-s} |h|^s   \int _{- \infty}^x |y|^s |u(y)|  \ \|   \mathfrak{m} _1^{- }(y,z ) -e_2\|   _{   L^2_z (\R )  } dy  \\& \le   2^{1-s} |h|^s e^{\| u \| _{L^1}} \| u \| _{L^{2,s}(\R)}
	\int _{- \infty}^x |y|^s
	\langle y \rangle ^{-s}   |u(y)| dy\le  2^{1-s} |h|^s e^{\| u \| _{L^1}} \| u \| _{L^{2,s}(\R)} \| u \| _{L^1}.\nonumber
\end{align}
Then, by  \eqref{eq:jost1}--\eqref{eq:jost3} we get  $\|    \mathfrak{m} _1^{- }(x,z+h)   -  \mathfrak{m} _1^{- }(x,z ) \|   _{   L^2_z (\R )  }
\le C  |h|^s  \| u \| _{L^{2,s}(\R)}$
for $x\le 0$, where $C$ is a fixed constant for $\| u \| _{L^{2,s}(\R)} \le \kappa _0$, for a preassigned bound $ \kappa _0$.
This implies that  for all $x\le 0$  we have    $\|    \mathfrak{m}_1^{- }(x,z )  -e_1   \|   _{   \dot H^s_z (\R )  }
\le C     \| u \| _{L^{2,s}(\R)}$ for some positive constant $C$.
Combined with \eqref{eq:jost-2}  this yields  the claim  \eqref{eq:jost0} for $j=1$ and for the minus sign.
The other cases are similar. \qed

\vspace{0.25cm}

Lemma \ref{lem:dirsc2}  provides   the direct scattering information we need.  Now we recall a number of facts
about inverse scattering.
The spectral data in the space
\begin{equation}\label{eq:spdata}
\mathcal{S} (s,n):= \left\{ r(z) \in H^s(\R), \quad (z_1,...,z_n) \in \mathbb{C}_+^n, \quad
(c_1,...,c_n) \in \C_*^n \right\}
\end{equation}
are used to recover the potential $u$ in matrix $Q(u)$ of the spectral system
(\ref{eq:zs}). Set
\begin{equation}\label{eq:vx}  \begin{aligned}  V_x(z):= \begin{pmatrix}
    1+|r(z)| ^2  &  e^{ -2 \im x z}\overline{r} (z) \\ e^{ 2 \im x z} {r} (z)
  & 1
 \end{pmatrix}
   \end{aligned}
\end{equation}
and consider the following Riemann--Hilbert (RH) problem:
\begin{itemize}
\item[(i)] $m(x,\cdot )$ is   meromorphic in $\mathbb{C}\backslash \R$;

\item[(ii)] $m(x,\cdot )$ has continuous boundary values $m_{\pm}(x,\cdot)$ on $\R$
satisfying
$$
m_+(x,z )=  m_-(x,z ) V_x(z);
$$

\item[(iii)] $\displaystyle \lim _{z\to \infty }m(x,z ) =1$;

\item[(iv)]  $m(x,z )$ has simple poles in $\mathcal{Z}=\mathcal{Z}_+\cup \mathcal{Z}_-$,
where $ \mathcal{Z}_-= \{ \bar{z}_1,...,\bar{z}_n \}$ in $\mathbb{C}_-$, and
for each $z_k\in  \mathcal{Z} _+ $ and $\bar{z}_k \in \mathcal{Z}_-$, we have
\begin{equation}\label{eq:resm}
\begin{aligned}
& \text{Res}_{z=z_k}m(x,z )= \lim _{z\to z_k}  m(x,z )
V_x(z_k), \\
& \text{Res}_{z=\overline{z}_k}m(x,z ) = \lim _{z\to \overline{z}_k}  m(x,z ) V_x(\overline{z}_k),
\end{aligned}
\end{equation}
with
\begin{equation}
\label{explicit-V-k}
V_x(z_k):= \begin{pmatrix}    0 &  0 \\  e ^{ 2\im x z_k } c_k   & 0 \end{pmatrix}, \quad
V_x(\overline{z}_k):=\begin{pmatrix}    0 &      -  e ^{ -2\im x \overline{z}_k }\overline{c  }_k\\ 0   & 0 \end{pmatrix}.
\end{equation}
\end{itemize}
>From the solution of the RH problem  (i)--(iv),  the potential $u$ in the matrix $Q(u)$
is found by means of {\it the reconstruction formula}:
\begin{equation}\label{eq:invscat1}
u(x) :=  2\im \lim _{z\to \infty }z  \, m _{12 }(x,z ).
\end{equation}

\begin{remark}
In terms of the analytic functions $\mathfrak{m}_{1,2}^{\pm }$ introduced
from the Volterra integral equations (\ref{eq:jost}), we have the equalities
$$
m_+(x,z) = \left[   {a(z)} ^{-1} \mathfrak{m}_{1}^-(x,z),   \mathfrak{m}_{2}^+(x,z) \right],
\quad m_-(x,z) = \left[ \mathfrak{m}_{1}^+(x,z),   {\overline{a (\overline{z})^{-1}}}     \mathfrak{m}_{2}^-(x,z) \right].
$$
\end{remark}
We introduce  now  the
  Cauchy operator $C_\R$ acting on
functions $h(z) \in L^2(\mathbb{R})$,
\begin{equation}\label{eq:cauchy1}
(C_\R h)(z) = \frac{1}{2\pi \im }
\int _{\R}\frac{h(\zeta)}{\zeta - z} d\zeta, \quad  z\in \C \backslash \R \, ,
\end{equation}
with the boundary values
\begin{equation*}
(C_\R ^{\pm} h)(z)=\lim _{\varepsilon \searrow 0} \frac{1}{2\pi \im }
\int _{\R}\frac{h(\zeta)}{\zeta -(z \pm \im \varepsilon ) } d\zeta, \quad z\in   \R.
\end{equation*}
The solution $m(x,z)$  of the  RH problem  (i)--(iv) is given by the following formula:
\begin{equation}\label{eq:syst2}
m(x,z) = 1- \sum _{\zeta \in \mathcal{Z}} \frac{M_x (\zeta ) V_x (\zeta )  }{\zeta -z}  +
\frac{1}{2\pi \im } \int _{\R}\frac{M_x(\zeta)(V_x(\zeta)-1)}{\zeta - z} d\zeta,
\end{equation}
where  $M_x(z)$ is defined for $z \in \mathbb{R} \cup \mathcal{Z}$ in the space
$M _{2\times 2}(\C )$ of complex $2 \times 2$ matrices  and satisfies system
\eqref{eq:syst1}--\eqref{eq:syst11} written below.

Lemma \ref{lem:zhou1}   below   implies that the map $\mathcal{G} _n\cap L^{2,s}(\R) \to \mathcal{S}(s,n)$,
which is defined by Lemmas \ref{lem:dirsc1} and \ref{lem:dirsc2},  is one-to-one.
This result is due to  Zhou \cite{Z},   but  we prove it for completeness, following the
argument in  Lemma 5.2 \cite{DP}.

\begin{lemma}
\label{lem:zhou1}
Let $r\in H^s (\R )$  with $s>1/2$. Then,
for any $x\in \R$ there exists and is unique a solution
$M_x : \R \cup \mathcal{Z}\to  M _{2\times 2}(\C )$
of the following system of integral and algebraic equations:
\begin{equation}
\label{eq:syst1}
M_x(z) = 1 -
\sum_{\zeta \in \mathcal{Z}} \frac{M_x (\zeta ) V_x (\zeta )  }{\zeta -z}
+ \lim _{\varepsilon \searrow 0} \frac{1}{2\pi \im }
\int _{\R}\frac{M_x(\zeta)(V_x(\zeta)-1)}{\zeta-(z - \im \varepsilon ) } d\zeta, \quad z\in \R
\end{equation}
and
\begin{equation}
M_x (z) = 1 -
\sum _{\zeta \in \mathcal{Z}\backslash \{ z \}} \frac{M_x(\zeta) V_x(\zeta)}{\zeta - z}
+ \frac{1}{2\pi \im } \int _{\R}\frac{M_x(\zeta)(V_x(\zeta)-1)}{\zeta- z} d\zeta,
\quad z \in \mathcal{Z} \label{eq:syst11}
\end{equation}
such that $(M_x(z) -1) \in L^{2}_z(\R )$.
\end{lemma}

\proof  For the operator $C _{V_x}$  defined by     $C _{V_x}h:=C^-(h(V_x-1))$,
the system of integral and algebraic equations (\ref{eq:syst1}) and (\ref{eq:syst11})
reduces to
\begin{equation}
\label{eq:syst20}
  (1- C _{V_x})  (M_x -1) +
\sum_{\zeta \in \mathcal{Z}} \frac{M_x (\zeta ) V_x (\zeta )  }{\zeta -z} = C _{V_x} 1, \quad z\in \R
\end{equation}
and
\begin{equation} \label{eq:syst21} \begin{aligned}&
M_x(z) +
\sum _{\zeta \in \mathcal{Z}\backslash \{ z \}} \frac{M_x (\zeta) V_x (\zeta)}{\zeta - z}
- \frac{1}{2\pi \im } \int _{\R}\frac{(M_x(\zeta)-1)(V_x(\zeta)-1)}{\zeta - z} d \zeta  \\
& \phantom{text} = 1 + \frac{1}{2\pi \im } \int _{\R}\frac{ (V_x(\zeta)-1)}{\zeta - z} d\zeta,
\quad z \in \mathcal{Z}.\end{aligned}
\end{equation}
By Lemma 5.2 \cite{DP}, there exists a fixed $c$ s.t.  for $\| r \| _{L^\infty (\R )}=\rho$   the operator
$ 1- C _{V_x} $  is invertible in $L^2 (\R )$  and $\|   (1- C _{V_x}) ^{-1}\| _{L^2\to L^2}\le c \langle \rho \rangle ^2.$
It is easy to conclude by the Fredholm alternative   that
the inhomogeneous system \eqref{eq:syst20}--\eqref{eq:syst21} admits exactly one solution
if  and only if  $f=0$ is the only  solution $f : \R \cup \mathcal{Z}\to  M _{2\times 2}(\C )$ with $f_{|\R }\in L ^{2}(\R ) $ of
\begin{equation}
  \label{eq:syhom}  \begin{aligned} &
	 (1- C _{V_x})  f +
\sum_{\zeta \in \mathcal{Z}} \frac{f (\zeta ) V_x (\zeta )  }{\zeta -z} = 0, \quad z\in \R \\&  f(z) +
\sum _{\zeta \in \mathcal{Z}\backslash \{ z \}} \frac{f (\zeta) V_x (\zeta)}{\zeta - z}
- \frac{1}{2\pi \im } \int _{\R}\frac{f(\zeta)(V_x(\zeta)-1)}{\zeta - z} d \zeta =0  ,
\quad z \in \mathcal{Z}. \end{aligned}
\end{equation}
We set  for $z\in \C  \backslash  (\mathcal{Z}\cup \R )$
\begin{equation}   \label{eq:syst3410}
m( z) := - \sum _{\zeta \in \mathcal{Z}} \frac{f (\zeta ) V_x (\zeta )  }{\zeta -z}  +
\frac{1}{2\pi \im } \int _{\R}\frac{f(\zeta)(V_x(\zeta)-1)}{\zeta -z} d\zeta.
\end{equation}
Notice that
\begin{equation}   \label{eq:syst341} \begin{aligned}
m_-( z)=-\sum _{\zeta \in \mathcal{Z}} \frac{f (\zeta ) V_x (\zeta )  }{\zeta -z}  +C^- (f (V_x -1)).
\end{aligned}
\end{equation}
By $C^- (f (V_x -1))=C_{V_x}f$ and \eqref{eq:syhom}
we get $m_-(z)=f(z)$  for $z\in \R $. We have
\begin{equation}   \label{eq:syst34} \begin{aligned}
m_+( z) &=  - \sum _{\zeta \in \mathcal{Z}} \frac{f (\zeta ) V_x (\zeta )  }{\zeta -z}  +C^+ (f (V_x -1)) \\
& =  - \sum _{\zeta \in \mathcal{Z}} \frac{f (\zeta ) V_x (\zeta )  }{\zeta -z}  +C^- (f (V_x -1))(z) +
f (z)(V_x(z) -1) \\
& =   m_-( z) + m_-( z) (V_x -1)   = m_-(z)V_x(z).
\end{aligned}
\end{equation}
We have
\begin{equation}  \label{eq:syst32} \begin{aligned}&
 0=\int _\R C^+(f (V_x-1))  \left (C^-(f (V_x-1))\right )^*  \\
 & =  \int _\R \left ( m_+( z) +  \sum _{\zeta \in \mathcal{Z}} \frac{f (\zeta ) V_x (\zeta )  }{\zeta -z} \right )   \left (m_-( z) +  \sum _{\zeta \in \mathcal{Z}} \frac{f (\zeta ) V_x (\zeta )  }{\zeta -z} \right )^* dz \\
 & = \int _\R  m_+ m_-^* dz + \sum _{\zeta \in \mathcal{Z}} \int _\R  \frac{ m_+( z) dz }{\overline{\zeta} -z}  V_x ^*(\zeta )f ^*(\zeta )
 + \left ( \sum _{\zeta \in \mathcal{Z}}  \int _\R \frac{ m_-( z) dz }{\overline{\zeta} -z}  V_x ^*(\zeta )f ^*(\zeta ) \right )^* \\& + \sum _{A\in \{ +,-\}}\, \sum _{\zeta , \xi  \in \mathcal{Z}_A}f (\zeta ) V_x (\zeta )V_x ^*(\xi )f ^*(\xi )\int _\R \frac {dz}{(\zeta -z)(\overline{\xi } -z)}       .
 \end{aligned}
\end{equation}
The $A=+$ term in the last line cancels with  the following:

\begin{equation*}   \begin{aligned}&    \sum _{\zeta \in \mathcal{Z}_+} \int _\R  \frac{ m_+( z) dz }{\overline{\zeta} -z}  V_x ^*(\zeta )f ^*(\zeta )   =
  \sum _{\zeta \in \mathcal{Z}_+} \int _\R  \frac{C^+ (f (V_x -1))(z)  dz }{\overline{\zeta} -z}  V_x ^*(\zeta )f ^*(\zeta ) -\\& - \sum _{\zeta \in \mathcal{Z}_+} \sum _{\xi  \in \mathcal{Z}} \int _\R  \frac{dz }{(\overline{\zeta } -z)  (\xi -z)}   f (\xi ) V_x (\xi )V_x ^*(\zeta )f ^*(\zeta )
\\& = - \sum _{\zeta ,\xi \in \mathcal{Z}_+} \int _\R  \frac{dz }{(\overline{\zeta } -z)  (\xi -z)}  f (\xi ) V_x (\xi )V_x ^*(\zeta )f ^*(\zeta ) .
 \end{aligned}
\end{equation*}
We have, by  $\text{Res}(m, \overline{\zeta})V_x ^*(\zeta ) =  f(\overline{\zeta}) V_x  (\overline{\zeta} )V_x ^*(\zeta )=-f(\overline{\zeta}) V_x ^*(\zeta )V_x ^*(\zeta )=0$,
   \begin{align}\nonumber &    \sum _{\zeta \in \mathcal{Z}_-} \int _\R  \frac{ m_+( z) dz }{\overline{\zeta} -z}  V_x ^*(\zeta )f ^*(\zeta ) =    \sum _{\zeta \in \mathcal{Z}_-} \int _\R  \left (    m_+( z) -\frac{\text{Res}(m, \overline{\zeta})}{z-\overline{\zeta}} \right ) \frac{V_x ^*(\zeta )f ^*(\zeta )    }{\overline{\zeta} -z}  dz
\\&\label{eq:syst35}= - 2\pi \im   \sum _{\zeta \in \mathcal{Z}_-}    f(\overline{\zeta})  V_x ^*(\zeta )f ^*(\zeta )  +   \sum _{\zeta \in \mathcal{Z}_-}   \sum _{ \xi  \in \mathcal{Z}_+ \backslash \{  \overline{\zeta}\}}  f(\xi )\underbrace{V_x  (\xi ) V_x ^*(\zeta )}_{0}f ^*(\zeta )\frac{2\pi \im}{\overline{\zeta} -\xi } \\&  = -2\pi \im   \sum _{\zeta \in \mathcal{Z}_-}    f(\overline{\zeta})  V_x ^*(\zeta )f ^*(\zeta ) .\nonumber
 \end{align}
Here we have used  the fact that for all $\zeta \in \mathcal{Z}$ by \eqref{eq:syhom}--\eqref{eq:syst3410}  we have
\begin{equation*}
\begin{aligned}&   \lim _{z\to \zeta } \left (    m ( z) -\frac{\text{Res}(m,  {\zeta})}{z- {\zeta}} \right )
= C( f(V_x-1))(\zeta ) - \sum _{\zeta' \in \mathcal{Z}\backslash \{ \zeta \}} \frac{f (\zeta ') V_x (\zeta')}{\zeta ' -\zeta }
=f(\zeta )   .
 \end{aligned}
\end{equation*}
The following term cancels with the $A=-$ term in  the last line of \eqref{eq:syst32}:
\begin{equation} \label{eq:syst36}  \begin{aligned}& \left (   \sum _{\zeta \in \mathcal{Z}_-} \int _\R  \frac{ m_-( z) dz }{\overline{\zeta} -z}  V_x ^*(\zeta )f ^*(\zeta ) \right )^*
=     \sum _{\zeta \in \mathcal{Z}_-} f  (\zeta ) V_x  (\zeta )[\int _{\R} \frac{(C^{-}(f (V_x-1)))^*}{\zeta -z} dz \\&  - \sum _{\xi  \in \mathcal{Z}}  \int _{\R} \frac{dz}{(\zeta -z) (\overline{\xi} -z)} f  (\zeta ) V_x  (\zeta )  V_x ^* (\xi ) f ^* (\xi ) ]
\\&
=  -  \sum _{\zeta , \xi  \in \mathcal{Z}_-}  f  (\zeta ) V_x  (\zeta )  V_x ^* (\xi ) f ^* (\xi ) \int _{\R} \frac{dz}{(\zeta -z) (\overline{\xi} -z)}  .
 \end{aligned}
\end{equation}
We have
   \begin{align}\nonumber & \left (   \sum _{\zeta \in \mathcal{Z}_+} \int _\R  \frac{ m_-( z) dz }{\overline{\zeta} -z}  V_x ^*(\zeta )f ^*(\zeta ) \right )^*=  \left (  \sum _{\zeta \in \mathcal{Z}_+} \int _\R  \left (    m_-( z) -\frac{\text{Res}(m, \overline{\zeta})}{z-\overline{\zeta}} \right ) \frac{V_x ^*(\zeta )f ^*(\zeta )    }{\overline{\zeta} -z}  dz\right )^*
\\&=      \left (   2\pi \im   \sum _{\zeta \in \mathcal{Z}_+}    f ( \overline{{\zeta}})  V_x ^* (\zeta )f ^*( {\zeta} ) \right )^* =       -2\pi \im   \sum _{\zeta \in \mathcal{Z}_+}    f ( {\zeta})  V_x  (\zeta )f ^*(\overline{\zeta} ) . \label{eq:syst37}
 \end{align}
The terms from  \eqref{eq:syst35} and  \eqref{eq:syst37} cancel out in \eqref{eq:syst32} because of $V_x  (\overline{\zeta} ) =-V_x ^* (\zeta )$.  Then,  by  $m_+ =m_-V_x$ , \eqref{eq:syst32}  yields
 \begin{equation}  \label{eq:syst38} \begin{aligned}&
 0=   \int _\R  m_-(z)V_x(z) m_-^* (z) dz   .
 \end{aligned}
\end{equation}
Since $V_x(z) $  is strictly positive, this implies $m_-(z)=f(z)=0$ for $z\in \R$.
  But then by  \eqref{eq:syst341}
we have also
$f (\zeta ) V_x (\zeta )=0$ for $\zeta \in \mathcal{Z}$. Then  $f (z )=0$ for $z \in \mathcal{Z}$ by \eqref{eq:syhom}.
So we have completed the proof that if $f$ solves  \eqref{eq:syhom} then $f=0$. \qed

\vspace{0.25cm}

We now recall another result  due to Zhou \cite{Z} on  the inverse scattering,
which we state in Lemma  \ref{lem:uc} below. This result is only stated for
the case of pure radiation solutions of the cubic NLS equation with $n = 0$.
We need  Lemma  \ref{lem:uc} in order to establish the fact that the map
$\mathcal{G} _0\cap  L^{2,s}(\R) \to \mathcal{S}(s,0)$  is not only
one-to-one but also onto.

\begin{lemma} \label{lem:uc}  Let   $r\in H^s(\R)$, $\mathcal{Z}=\emptyset $,
and consider the potential $u$ defined by the reconstruction formula \eqref{eq:invscat1}.
Then $u \in L^{2,s}(\R)$. Furthermore,  for any positive $\kappa _0$, there is a constant $C$
such that for $\| r\|_{L^\infty (\R)} \le  \kappa _0$,
we have $\| u\| _{ L^{2,s}(\R)}\le C\| r\| _{H^s(\R)} $.
\end{lemma}

\proof We only sketch the argument, referring to references \cite{Z,DZ,DP}
for more information  and details. We first sketch $u(x) \in L^{2,s}(\R _+)$.
We factorize   the matrix in  \eqref{eq:vx}  writing  $V_x(z) =  V _{x-}  ^{-1}  V_{x+} $, where
\begin{equation}\label{eq:vxpm}  \begin{aligned}
V_{x+}(z):=\begin{pmatrix}
    1   & 0 \\  e^{ 2\im xz } {r} (z)
  & 1
 \end{pmatrix} \, , \,   V _{x-} (z):= \begin{pmatrix}
    1   & -e^{  -2\im xz }   \overline{r} (z) \\  0
  & 1
 \end{pmatrix}
  . \end{aligned}
\end{equation}
Set  now  $C _{w_x}h:=C^+(hw _{x - })+ C^-(hw _{x + }) $ for
$w _{x \pm }  := \pm (V _{x \pm } -1)$. Then we consider a function $\mu _x\in 1+L^2(\R)$ such that
\begin{align} \label{eq:syst11a} &
 (1-C _{w_x}) (\mu _x)   (z)= 1 .
 \end{align}
For $w_x(\zeta ):= V _{x+}(\zeta )-V _{x-}(\zeta )$ we get that
the $m(x,z)$ in \eqref{eq:syst2} (in the case when all the $c_j=0$)
can be expressed also as
\begin{equation}\label{eq:syst2a}\begin{aligned} &
m(x,z)=   1+  \frac{1}{2\pi \im } \int _{\R}\frac{\mu _x(\zeta )w_x(\zeta )}{\zeta -z} d\zeta .
 \end{aligned}
\end{equation}
For $x\ge  0$,  by the argument in Lemma 3.4
\cite{DZ}   for a fixed $c_s$    we have
\begin{equation*}  \begin{aligned} &  \|  C^\pm (1- V _{x\mp} (z) )\| _{L^2_z(\R )}\le
c_s \langle x \rangle ^{-s} \|  r\| _{H^s(\R)}
\end{aligned}\end{equation*}
(notice that   $x\le  0$  in Lemma 3.4 \cite{DZ}, because of
the different definition of the operator in \eqref{eq:zs}).
This implies immediately
$$
\| C _{w_x}1\| _{L^2_z(\R )}\le 2c_s \langle x \rangle ^{-s} \|  r\| _{H^s(\R)}.
$$
We   consider
\begin{equation*}  \begin{aligned} &     \mu _x  -1= \left ( 1-C _{w_x}  \right ) ^{-1}       C _{w_x}1
\end{aligned}\end{equation*}
and correspondingly
\begin{equation*}  \begin{aligned} &   \|  \mu _x  -1\| _{L^2_z}   \le
\| \left ( 1-C _{w_x}  \right ) ^{-1} \| _{L^2_z \to L^2_z}     \|  C _{w_x}1 \| _{L^2_z}  .
\end{aligned}\end{equation*}
We have $\| \left ( 1-C _{w_x}  \right ) ^{-1} \| _{L^2_z \to L^2_z} \le c\langle \rho \rangle ^{2}$
by Lemma 5.2 \cite{DP} for a fixed $c$, where $\rho := \| r\| _{L^\infty(\R)}$.
We conclude that for $x\ge  0$ and for any $\kappa _0$ there is a constant $C$
such that
\begin{equation*}  \begin{aligned} &   \| \mu _x  -1\| _{L^2_z}   \le C \langle x \rangle ^{-s} \|  r\| _{H^s(\R)}
\end{aligned}\end{equation*}
for $\rho \le  \kappa _0$. Finally, the argument in Theorem 3.5
\cite{DZ}   yields $\|   u    \| _{L^{2,s}(\R _+)}\le C \|  r\| _{H^s(\R)}$.

In order to prove   $u(x) \in L^{2,s}(\R _-)$   we consider instead the decomposition
\begin{equation*}  \begin{aligned} &  V_x= \begin{pmatrix}
    1   & 0 \\   \frac{e^{   2\im  xz  } {r} (z)}{1+|{r} (z)|^2}
  & 1
 \end{pmatrix}    \begin{pmatrix}
   1+|{r} (z)|^2   & 0 \\   0
  &  \frac{1}{1+|{r} (z)|^2}
 \end{pmatrix} \begin{pmatrix}
    1   & \frac{e^{  -2\im  xz }   \overline{r} (z)}{1+|{r} (z)|^2} \\  0
  & 1
 \end{pmatrix}.
\end{aligned}\end{equation*}
We then consider the RH problem with matrix
$\widetilde{V}_x:= \delta _-^{\sigma _3} {V}_x \delta _ +^{-\sigma _3}$  for $\delta (z)$ the solution
of the problem \eqref{eq:scalrRH} with $z_0=+\infty $  introduced later in
Proposition \ref{prop:scalarRH}. Correspondingly we get estimates
$\|    \widetilde{{u}}    \| _{L^{2,s}(\R _-)}\le C \|  \widetilde{r}\| _{H^s(\R) }\le c \|   {r}\| _{H^s(\R) }$
for a function $\widetilde{{u}}$  associated to $\widetilde{r} :=r \delta _+\delta _- $  and
for fixed $c$   when $\rho \le \kappa _0$, by proceeding as above. Finally, $\widetilde{u} = {u} $.
For more details see \cite{DZ}. \qed

\vspace{0.25cm}

We now discuss  the representation of the solutions of the Cauchy problem  \eqref{eq:nls}
in terms of the inverse scattering transform. We  recall the following result, see \cite{GV}.

\begin{theorem}
\label{th:GV} Given $u_0\in H^1(\R ) $,  then there exists a unique solution
$u(t) \in C^0 (\R, H^1(\R )) \cap  L^4_{loc}(\R, L^\infty(\R ))$
of the integral equation (\ref{eq:nls1}).
 \end{theorem}

The solution of Theorem  \ref{th:GV} is the same of   Theorem  \ref{th:L2}.

Suppose   that $u_0 \in   H^1(\R )\cap L^{2,s}(\R)$  for fixed $s \in \left(\frac{1}{2},1\right]$.
Then the solution remains in  $u(t)\in    H^1(\R )\cap L^{2,s}(\R) $   for all $t \in \R$,
by standard arguments (see p. 1072 in \cite{MS}, which can be extended to non-integer $s$
by Lemma 2.3 in \cite{HN}). For the solution of the cubic NLS equation (\ref{eq:nls})
with $u_0 \in H^1(\R) \cap L^{2,s}(\R)$,  the time evolution of the scattering data is well-defined,
according to the following result (see p.39 in \cite{ABT}):

\begin{lemma} \label{lem:invsc}
For an  initial datum  $u_0\in H^1(\R )\cap  L^{2,s}(\R) \cap \mathcal{G}$
we have $u(t)\in H^1(\R )\cap  L^{2,s}(\R) \cap \mathcal{G}$   for all $t \in \R$
and the  spectral data $\mathcal{S}(s,n)$ in  \eqref{eq:spdata} evolve  as follows:
\begin{equation}\label{eq:spdatau}
e^{ 4\im z^2 t} r(  z ) \in H^s(\R), \quad (z_1,...,z_n) \in \mathbb{C}_+^n, \quad
(e^{ 4\im z_1^2 t}c_1,...,e^{ 4\im z_n^2 t}c_n) \in \C_*^n.
\end{equation}
\end{lemma}


\begin{remark}
\label{lem:soliton}
To recover the solitons (\ref{eq:groundstate}), we take the spectral data:
\begin{equation}\label{eq:soliton1}
r=0, \quad  z_1 = \alpha _1 + \im \beta _1 \in \mathbb{C}_+, \quad  e^{-4\im z_1^2 t}c_1 \in \C_*.
\end{equation}
Then, we obtain
\begin{equation}\label{eq:soliton2}
u(x,t) = -2\im \beta _1  e^{-2\im \alpha _1x-4\im t (\alpha^2_1- \beta ^2_1)-\im  \psi _0}
 \text{sech}(2\beta _1 x+8t\alpha _1\beta _1- \delta _0 ),
\end{equation}
where $\delta_0 :=  \log \left (\frac{|c_1|}{ 2\beta _1}\right )$ and $\psi _0:=\arg  (c_1)$.
Note the correspondence: $\omega = 2 \beta_1$ and $v = -2 \alpha_1$,
for solitons in (\ref{eq:groundstate}).
\end{remark}

\section{Dispersion for pure radiation solutions}
\label{sec:purerad}

Elements of  ${\mathcal G}$ such that ${\mathcal Z}_+=\emptyset $ generate pure radiation
solutions of the cubic NLS equation. These solutions satisfy the following asymptotic behavior.

\begin{theorem}
\label{th:purerad} Fix $s\in (1/2, 1]$.
Let   $u_0\in  {\mathcal G}\cap L^{2,s}(\R )$ such that
  ${\mathcal Z} =\emptyset $.
Then  there exist   constants $C(u_0 )>0$  and    $T(u_0 )>0$   such that
the solution of the cubic NLS equation \eqref{eq:nls}    satisfies
\begin{equation*} \begin{aligned} &
\|    u(t , \cdot )  \| _{L^\infty ( \R)} \le
C(u_0) |t|^{-\frac 12}  \text{ for all $|t|\ge T(u_0)$.}
   \end{aligned}
\end{equation*}
There are furthermore constants $C_0 > 0$, $T_0>0$ and small $\varepsilon_0 > 0$ such that
for   $ \|    u_0 \| _{L^{2,s}(\R)}<\varepsilon_0 $, we can take
$C(u_0) = C_0 \|    u_0 \| _{L^{2,s}(\R)}$   and $T(u_0)=T_0$.
\end{theorem}

\begin{remark}
In \cite{HN}, the result of Theorem \ref{th:purerad} is proved with $L^{2,s}(\R)$
replaced by $\Sigma _s$ for any $s> \frac{1}{2}$,  only in the case
of small $u_0$ with  $  \|    u_0 \| _{\Sigma _s } < \varepsilon _0 $.
In the case of the defocusing NLS equation \eqref{eq:nls} (that is, with $+ 2  |u|^2  u$  replaced by
$- 2  |u|^2  u$), Theorem \ref{th:purerad} for $s=1$
is proved in  \cite{DZ2,DZ}. For the focusing NLS equation \eqref{eq:nls},
Theorem \ref{th:purerad} for $s=1$ is proved in \cite{DP}. Notice also that
all these references contain  proofs of the asymptotic expansions  for
the solution $u$ at large $t$, which we do not discuss here.
\end{remark}

In the  rest of Sect. \ref{sec:purerad} we prove  Theorem \ref{th:purerad}. With minor modifications,  we
follow  closely    the  proof in  \cite{DM},
which  involves  the $\overline{\partial}$  operator, where
$
\overline{\partial} :=  \frac 12 (\partial _x +\im \partial _y).
$
Here we extend the result in \cite{DM}, valid for $s=1$, to any $s\in (1/2, 1]$.

\subsection{Proof of {Theorem}
\ref{th:purerad}}
\label{ssec:purerad}

The proof starts by assuming additionally that  $u_0\in H^1 (\R )$.
Fix   $z_0\in \R .$
First of all   we consider the  scalar RH problem

\begin{equation} \label{eq:scalrRH}
\left\{\begin{matrix}
          \delta _+ (z)  = \delta _- (z)  (1+|{r} (z)|^2)  \text{ for $z<z_0$} \,\\
  \delta _+ (z)  = \delta _- (z)     \text{ for $z>z_0$}
\end{matrix}\right.
\end{equation}
with $\delta (z)$ holomorphic in $\C \backslash \R$ and
$\delta (z)\to 1$ as $z\to \infty $.  The following statement is in Prop. 5.1 \cite{DP}  and is elementary to prove.
\begin{proposition}\label{prop:scalarRH}
 We have
\begin{equation*}  \begin{aligned}      \delta  (z)   =    e^{  \gamma (z)}\, , \,  \gamma (z):= { \frac 1{2\pi \im } \int _{-\infty}^{z_0}    \frac{\log (1+|{r} (\varsigma )|^2)}{\varsigma -z} d\varsigma}.
   \end{aligned}
\end{equation*}
For $z\not \le z_0$ we have $  \delta  (z)  = \overline{\delta  (\overline{z})}$ and    $\langle \rho \rangle ^{-1}
\le   | \delta   (z)  |   \le \langle \rho \rangle   $  where $\rho :=\| r \| _{L^\infty(\R)}$;   for $\mp \Im z >0$ we have   $ | \delta ^{\pm 1} (z)  |\le 1$.

\end{proposition}

The function $ \gamma (z)$ has an expansion
\begin{equation} \label{eq:gamma}  \begin{aligned}     \gamma (z)  =    \im \nu (z_0) \log (z-z_0) +   \im \nu (z_0)(z-z_0) \log (z-z_0) \\-\im \nu (z_0)(z-z_0+1) \log (z-z_0+1) + \beta (z, z_0)
   \end{aligned}
\end{equation}
where in the r.h.s. the main term is the first, and where
\begin{equation} \label{eq:beta}  \begin{aligned}   &  \nu (z_0)  := -\frac{1}{2\pi}   \log (1+|r(z_0)|^2)   \text{ and for }\chi  (\zeta,z_0)= \chi _{[z_0-1,z_0]}(\zeta ) (\zeta -z_0+1) \\&    \beta (z, z_0) =  \int _{-\infty}^{z_0}  \left \{  {\log (1+|{r} (\zeta)|^2)} -   \log (1+|r(z_0)|^2) \chi  (\zeta,z_0)    \right \}  \frac{d\zeta}{2\pi \im (\zeta -z)} .
   \end{aligned}
\end{equation}

Let
\begin{equation}\label{eq:theta}  \begin{aligned}    \text{$\theta (z):= 2(z-z_0)^2-2z_0^2$ with $z_0:=- \frac x {4t}.$}
   \end{aligned}
\end{equation}
Then we consider the RH problem  (i)--(iii) with
\begin{equation*}  \begin{aligned}   v_{t,x}(z):= \begin{pmatrix}
    1+|r(z)| ^2  & e^{ -2\im t \theta }   \overline{r} (z) \\  e^{ 2\im t \theta } {r} (z)
  & 1
 \end{pmatrix}
  . \end{aligned}
\end{equation*}
We factorize

\begin{equation*}    v_{t,x}  =
\left\{\begin{matrix}
   W_LW_R \text{ for $z>z_0$} \,\\
 U_LU_0U_R \text{ for $z<z_0$}
\end{matrix}\right.
\end{equation*}
 for
\begin{equation} \label{eq:jump1} \begin{aligned}  &  W_L=  \begin{pmatrix}
    1   & e^{ -2\im t \theta }   \overline{r} (z) \\  0
  & 1
 \end{pmatrix}    \, , \,    W_R =  \begin{pmatrix}
    1   & 0 \\  e^{2\im t \theta } {r} (z)
  & 1
 \end{pmatrix}     \, , \,  U_L=  \begin{pmatrix}
    1   & 0 \\  \frac{e^{ 2\im t \theta } {r} (z)}{1+|{r} (z)|^2}
  & 1
 \end{pmatrix}  \\&   U_R = \begin{pmatrix}
    1   & \frac{e^{ -2\im t \theta }   \overline{r} (z)}{1+|{r} (z)|^2} \\  0
  & 1
 \end{pmatrix}  \, , \,      U_0 = \begin{pmatrix}
   1+|{r} (z)|^2   & 0 \\   0
  &  \frac{1}{1+|{r} (z)|^2}
 \end{pmatrix} . \end{aligned}
\end{equation}

We end   Sect. \ref{ssec:purerad} with an estimate on the function $\beta (z,z_0)$.

\begin{lemma}
  \label{lem:beta} Let $L_\phi = z_0+e^{-\im \phi}\R   =\{ z=z_0+e^{-\im \phi}u:u\in \R      \}$.   Consider the $s\in (1/2, 1]$ in Theorem \ref{th:asstab}.
  Then      there is a fixed $C(\rho ,s)$ s.t.  for any $z_0\in \R$ and any $\phi \in (0,\pi )$

    \begin{align} \label{eq:beta1}   &   \| \beta ( e^{-\im \phi} \cdot , z_0) \| _{H^s(\R )} \le C(\rho ,s)   \| r \| _{H^s(\R )} \\& \label{eq:beta2}   | \beta ( z , z_0)- \beta ( z_0 , z_0)  | \le C(\rho ,s)     \| r \| _{H^s(\R )} |z-z_0|^{s-\frac 12}  \text{ for all $z\in L_\phi$}.
   \end{align}

 \end{lemma}
\proof  First of all  these estimates hold for $s=1$, and are a consequence of  $\| C_\R f\| _{H^\tau (L_\phi )} \le C_\tau  \|  f\| _{H^\tau (\R )}  $ for $\tau =0,1$, which are
proved in  Lemma 23.3   \cite{BDT}.  We obtain \eqref{eq:beta1} for $\tau =s$
when $s\in  (0, 1)$  by interpolation.  The estimate    \eqref{eq:beta2}
is   a consequence of  \eqref{eq:beta1} and of the following elementary estimate when  $s\in (1/2, 1]$:
 \begin{align} \label{eq:beta3}   | f (x)- f ( y)  | \le C_s     \| f \| _{H^s(\R )} |x-y|^{s-\frac 12}  \text{for all $x,y\in \R$ and $f\in H^s(\R )$  for a fixed $C_s$}.
   \end{align}
  This is an elementary consequence of $ f (x+h)- f ( x) =\frac{1}{\sqrt{2\pi}} \int e^{\im x \xi}
  (e^{\im h \xi}-1) \widehat{f}(\xi ) d\xi $  for $y=x+h$. Then for any $\kappa >0$ we have
  for a fixed $C_s$ \begin{equation*}  \begin{aligned} &
 |f (x+h)- f ( x)|\le  \frac{|h|}{\sqrt{2\pi}}  (\int _{|\xi |\le \kappa }  |\xi | ^{2-2s} d\xi ) ^{\frac{1}{2}} \| f \| _{H^s} +   \frac{1}{\sqrt{2\pi}}  (\int _{|\xi |\ge \kappa }  |\xi | ^{ -2s} d\xi ) ^{\frac{1}{2}} \| f \| _{H^s} \\& \le C _s (|h|  \kappa ^{\frac{3-2s}{2}}
 + \kappa ^{\frac{1-2s}{2}}) \| f \| _{H^s}.
\end{aligned}
\end{equation*}
   The r.h.s. equals  $2C_s|h|^{s-\frac 12} \| f \| _{H^s}$ for $\kappa  =|h|^{ -1} $.

\qed

\subsubsection{The  model RH problem}
\label{ssec:model}
We   consider the  RH problem

\begin{equation} \label{eq:model1}
\left\{\begin{matrix}    P \text{ analytic in } \C \backslash \Sigma _P \\
    P(\zeta )=  1+ \frac{P_1}{\zeta}  +  O(\zeta ^{-2})  \text{ as }   \zeta \to  \infty \\ P_+(\zeta )=P_+(\zeta ) V_P(\zeta ) \text{  in } \Sigma _P
\end{matrix}\right.
\end{equation}
where    $\Sigma _P:= \overline{\cup _{n=1}^{4}\Sigma _P^{n}}$  with
$ \Sigma _P^{1}=  e^{ \im  \frac{\pi}{4}  } \R_+$,   $ \Sigma _P^{2}=  e^{- \im  \frac{\pi}{4}  } \R_-$, $ \Sigma _P^{3}=  e^{ \im  \frac{\pi}{4}  } \R_-$  and    $ \Sigma _P^{4}=  e^{ -\im  \frac{\pi}{4}  } \R_+$ inheriting the orientations of $\R _\pm $.
 The   matrix $V_P(\zeta )$
is defined by

\begin{equation}\label{eq:defVP} V_P(\zeta ):=
\left\{\begin{matrix}
\begin{pmatrix}
    1   & 0 \\       r_0\zeta ^{-2\im \nu _0} e^{\im \zeta ^2/2}
  & 1
 \end{pmatrix}  \text{ for }   \zeta \in  \Sigma _P^{1}   \\
 \begin{pmatrix}
    1   & \frac{ \overline{r } _0}{1+|r_0|^2}\zeta ^{ 2\im \nu _0} e^{-\im \zeta ^2/2} \\      0
  & 1
 \end{pmatrix}  \text{ for }   \zeta \in  \Sigma _P^{2}    \\\begin{pmatrix}
    1   & 0 \\     \frac{r_0}{1+|r_0)|^2}\zeta ^{ -2\im \nu _0} e^{ \im \zeta ^2/2}
  & 1
 \end{pmatrix}  \text{ for }   \zeta \in  \Sigma _P^{3}   \\
     \begin{pmatrix}
    1   &  \overline{r } _0 \zeta ^{ 2\im \nu _0} e^{-\im \zeta ^2/2} \\      0
  & 1
 \end{pmatrix}  \text{ for }   \zeta \in  \Sigma _P^{4}
\end{matrix}\right.
\end{equation}
where $r_0  $   is a free parameter
that  we fix in \eqref{eq:defVP1}  
and  $\nu _0=\nu (z_0)$.
The solution of this RH problem can be  worked out following word by word  \cite{DZ2} pp.54--57.
Set
\begin{equation}   \label{eq:kdef} \begin{aligned} &
   k_1:=
    \frac{-\im \sqrt{2\pi} e^{\im \pi /4} e^{-\pi \nu _0 /2}}{r_0\Gamma (-\im \nu _0)} \, , \,  k_2:=\frac{\nu _0}{k_1}  .
\end{aligned} \end{equation}
Consider for $\Im \zeta > 0$  the matrix $\Psi  ^+ (\zeta )$ with
\begin{equation*} \begin{aligned} &\Psi  ^+ _{11}    (\zeta )    =
  e  ^{-3 \pi  \nu   _0/4 } D_{\im \nu _0} (e  ^{-3 \im \pi    /4 } \zeta  )  \, , \quad
	\Psi  ^+ _{22}    (\zeta )    = e  ^{  \pi  \nu    _0/4 } D_{-\im \nu _0 } (e  ^{-  \im \pi    /4 } \zeta  )  \ , \\&
	\Psi  ^+ _{12}    (\zeta )    =       e  ^{  \pi  \nu   _0 /4 }  (-\im k_2) ^{-1}
	\left (  \partial _\zeta (  D_{-\im \nu  _0 }     (e  ^{-  \im \pi    /4 } \zeta  ) )   -    \frac{\im \zeta }2  D_{-\im \nu _0 }     (e  ^{-  \im \pi    /4 } \zeta  )    \right )  \ , \\
  &  	\Psi  ^+ _{21}    (\zeta )    =       e  ^{-3 \pi  \nu   _0 /4 }  (\im k_1) ^{-1}
	\left (  \partial _\zeta (  D_{ \im \nu _0 }     (e  ^{-  \im 3\pi    /4 } \zeta  ) )   +    \frac{\im \zeta }2  D_{ \im \nu _0 }     (e  ^{-  \im 3\pi    /4 } \zeta  )    \right ) \, .
 \end{aligned}
\end{equation*}
Consider for $\Im \zeta  < 0$  the matrix $\Psi  ^- (\zeta )$ with
\begin{equation*} \begin{aligned} &\Psi  ^- _{11}    (\zeta )    = e  ^{  \pi  \nu  _0/4 } D_{ \im \nu  _0 } (e  ^{   \im \pi   /4 } \zeta  ) \, , \quad
	\Psi  ^- _{22}    (\zeta )    =
 e  ^{-3 \pi  \nu    _0 /4 } D_{\im \nu   _0} (e  ^{ 3 \im \pi   /4 } \zeta  )  \, ,
\\&
	\Psi  ^- _{12}    (\zeta )    =       e  ^{-3 \pi  \nu    _0 /4 }   (-\im k_2) ^{-1}  \left (  \partial _\zeta (  D_{- \im \nu   _0}     (e  ^{  \im 3\pi   /4 } \zeta  ) )   -    \frac{\im \zeta }2  D_{- \im \nu   _0}     (e  ^{  \im 3\pi  /4 } \zeta  )    \right ) \, ,  \\
	&  	\Psi  ^- _{21}    (\zeta )    =    e  ^{  \pi  \nu /4 }  (\im k_1) ^{-1}
	\left (  \partial _\zeta (  D_{ \im \nu   _0}     (e  ^{   \im \pi  /4 } \zeta  ) )   +   \frac{\im \zeta }2  D_{ \im \nu   _0}     (e  ^{   \im \pi   /4 } \zeta  )    \right ) \, .
 \end{aligned}
\end{equation*}
Here   $ D_{ a} (\zeta )  $  is the  unique  entire function       solving
\begin{equation*}   \begin{aligned} &   \frac {d^2}{d \zeta ^2}D_{ a} (\zeta )   +
\left (\frac 12 -\frac{\zeta ^2} 4+a  \right ) D_{ a} (\zeta ) =0
\end{aligned} \end{equation*}
that for $|\arg (\zeta ) |<3\pi /4$   satisfies
\begin{align} &  \nonumber  D_{ a} (\zeta )   \sim  e^{   -\zeta ^2/4 } \zeta ^a   \left (  1 +\sum _{n=1} ^\infty (-1)^n 2^{ n}  \frac{ \prod _{j=1}  ^n \left (  \frac a2 -  (j-1) \right ) \left (  \frac a2 -  (j-1/2) \right )   } {   n!   \zeta   ^{2n} }   \right ) \\&\text{and such that }   D_{ a} ' (\zeta ) +\frac{\zeta}{2} D_{ a} (\zeta ) =aD_{ a-1} (\zeta ),\label{eq:idD}
\end{align}
see Chapter 16  \cite{WW}.
If we introduce the angular sectors
\begin{equation*}   \begin{aligned} & \Omega _1= \{ \zeta : \arg  \zeta \in (0,    \pi /4 \} \, , \quad  \Omega _2= \{ \zeta : \arg  \zeta \in ( \pi /4, 3   \pi /4) \}  \ , \\&
\Omega _3= \{ \zeta : \arg  \zeta \in (3 \pi /4,  \pi  ) \} \, , \quad  \Omega _4= \{ \zeta : \arg  \zeta \in (\pi ,  5\pi /4 ) \}  \ , \\&\Omega _5= \{ \zeta : \arg  \zeta \in (\pi +   5\pi /4,  7\pi /4) \} \, ,\quad \Omega _6= \{ \zeta : \arg  \zeta \in (7\pi /4 , 2\pi  ) \} \  ,
\end{aligned} \end{equation*}
then, following \cite{DZ2}, see also \cite{DP}, we have
\begin{equation}  \label{eq:defP} \begin{aligned} &
P(\zeta ) =\Psi  ^+ (\zeta )     \begin{pmatrix}
    1   & 0 \\      - {r}  _0
  & 1
 \end{pmatrix} \zeta ^{ -\im \nu   _0\sigma _3} e^{ \im \zeta ^2\sigma _3/4}  \text{ for } \zeta \in \Omega _1 \, ,\\&
P(\zeta ) =\Psi  ^+ (\zeta )  \zeta ^{ -\im \nu   _0\sigma _3}e^{ \im \zeta ^2\sigma _3/4}
 \text{ for } \zeta \in \Omega _2 \, ,\\&
P(\zeta ) =\Psi  ^+ (\zeta )     \begin{pmatrix}
    1   & \frac{ -\overline{{r}} _0}{1+|{r}  _0|^2}  \\      0
  & 1
 \end{pmatrix} \zeta ^{ -\im \nu   _0\sigma _3}e^{ \im \zeta ^2\sigma _3/4}
 \text{ for } \zeta \in \Omega _3 \, ,\\&
P(\zeta ) =\Psi  ^- (\zeta )     \begin{pmatrix}
    1   &0  \\       \frac{  {{r}}  _0}{1+|{r}  _0|^2}
  & 1
 \end{pmatrix} \zeta ^{ -\im \nu   _0\sigma _3}e^{ \im \zeta ^2\sigma _3/4}
 \text{ for } \zeta \in \Omega _4 \, ,\\&
P(\zeta ) =\Psi  ^- (\zeta )  \zeta ^{ -\im  \nu   _0\sigma _3}e^{ \im \zeta ^2\sigma _3/4}
 \text{ for } \zeta \in \Omega _5 \, ,\\&
P(\zeta ) =\Psi  ^- (\zeta )   \begin{pmatrix}
    1   &  { \overline{{r}} _0}   \\      0
  & 1
 \end{pmatrix}  \zeta ^{ -\im \nu   _0\sigma _3} e^{ \im \zeta ^2\sigma _3/4}
 \text{ for } \zeta \in \Omega _6 \ .
\end{aligned} \end{equation}
The fact that $P(\zeta )$ satisfies the model RH problem \eqref{eq:model1}
can  be seen by direct computation  (specifically, it solves an equivalent  RH with
an additional  jump matrix 1  over $\R$: this fact can be checked directly  by exploiting
the fact that $\Psi  ^+ (\zeta )=\Psi  ^- (\zeta ) \begin{pmatrix}
    1+|{r}_0| ^2  &   \overline{r}_0 \\ {r}_0
  & 1
 \end{pmatrix}$ and the monodromy properties of $z ^{\nu}$ like in p. 48 \cite{DZ2}).

By elementary computations which use \eqref{eq:defP} and   \eqref{eq:idD},   see \cite{DZ2},   we have
\begin{equation} \label{eq:p1}   \begin{aligned} &
 \lim _{\C _\pm \ni \zeta \to \infty }[\Psi  ^\pm (\zeta )  \zeta ^{ -\im  \nu   _0\sigma _3}e^{ \im \zeta ^2\sigma _3/4}-1] \zeta= P_1  \text{ with } P_1:= \begin{pmatrix}
    0   &  k_1  \\       k_2
  & 0
 \end{pmatrix} .
\end{aligned}
\end{equation}
Exploiting the rapid convergence to $1$ as $\zeta \to \infty$ of the extension of $(V_P)| _{\Sigma _P^1}$ to $\Omega _1$, of
$(V_P)| _{\Sigma _P^2}$ to $\Omega _3$, of $(V_P)| _{\Sigma _P^3}$ to $\Omega _4$ and of $(V_P)| _{\Sigma _P^4}$ to $\Omega _6$, it is easy to conclude that
$\displaystyle \lim _{\zeta \to \infty } \zeta (P(\zeta )-1) =
P_1  $ in each sector $\Omega _j$.
In each sector we have $\det P(\zeta )=1$, see p.54 \cite{DZ2}.

With respect to the analysis in \cite{DZ2} we need to add few
more remarks of quantitative nature on $P(\zeta )$.
We fix 
\begin{equation}\label{eq:defVP1} r_0:= \widehat{r}_0 e^{\im \nu _0 \log ( 8t) -4 \im t z_0^2} \text{ and } \widehat{r} _0:= r(z_0) e^{-2\im \nu (z_0) -2\beta (z_0, z_0)}.
\end{equation}
By    $|r _0|= |r (z_0)| \le C_s \| r\| _{ {H}^s(\R)} \le C (u_0 )$
there is a  $ C(u_0)$  such that by  \eqref{eq:beta} and \eqref{eq:kdef}  we get $|k_1|+|k_2|\le  C (u_0)$. Furthermore the following is true.

\begin{lemma} \label{lem:asP} Let $\rho = \| r\| _{L^\infty(\R)}$.  For any $\rho _0$  there exists a $C$ such that for  $\rho \le \rho _0$ we have
\begin{align} \label{eq:p2}&
| P(\zeta )  |\le  C \text{ for all $ \zeta \not \in \R$  and }  \\&
| P(\zeta )-1- P_1/\zeta |\le  {C \, \rho}{|\zeta |  ^{-2}   }\,    \text{ if also  $|\zeta |\ge 1$ } . \label{eq:p21}
\end{align}
\end{lemma}
\proof   We focus only on   \eqref{eq:p21}, since \eqref{eq:p2}  follows by \eqref{eq:p21} and
by the fact that $ D_{ a} (\zeta )  $  is an entire function in  $(a,\zeta )$.
The proof of  \eqref{eq:p21}   is based on formulas for $D_{\im \nu _0}(\zeta)$ for which we refer to
Chapter 16 \cite{WW}.

 Recall that  $D_{\im \nu _0}(\zeta)= 2^{\frac{{\im \nu _0}}{2}+\frac{1}{4}} \zeta ^{-\frac{1}{2}} W_{\frac{{\im \nu _0}}{2}+\frac{1}{4},-\frac{1}{4}} (\frac{\zeta ^2}{2})$, where for $|\arg (z)|<3\pi /2$   we have
\begin{equation*}  \begin{aligned} &  W_{\frac{{\im \nu _0}}{2}+\frac{1}{4},-\frac{1}{4}} (z) =  e^{-z /2}  z ^{\im \frac{{  \nu _0}}{2}+\frac{1}{4}}  \big  [ 1  -z^{-1}  \frac{\Gamma \left (\frac{3}{2}-\im \frac{\nu _0}{2}  \right )\Gamma \left  (1-\im \frac{\nu _0}{2}  \right )}{\Gamma \left (\frac{1}{2}-\im \frac{\nu _0}{2}  \right )\Gamma \left ( -\im \frac{\nu _0}{2}  \right )} +\\&   \frac 1{\Gamma \left (\frac{1}{2}-\im \frac{\nu _0}{2}  \right )\Gamma \left ( -\im \frac{\nu _0}{2}  \right )} \frac{1}{2\pi \im } \int _{-\im \infty -\frac{3}{2}}^{+\im \infty -\frac{3}{2}} z ^{\varsigma }{\Gamma \left (\varsigma   \right )\Gamma \left (-\varsigma +\frac{1}{2}-\im \frac{\nu _0}{2}  \right )\Gamma \left  (-\varsigma -\im \frac{\nu _0}{2}  \right )} d\varsigma  \big   ] . \end{aligned}
\end{equation*}
To bound  the integral we use:
\begin{equation*}  \begin{aligned} &   |z^{\varsigma }|= |z|^{\Re (\varsigma )}   e^{-t\arg (z)}
\text{ for
  $\varsigma = \Re (\varsigma ) +\im t$;}  \\&  |\Gamma (z)|\le   \sqrt{2\pi } |z^{z-\frac{1}{2}}|  e^{\frac{K}{  \Re (z )}}  \text{ for
  $\Re z>0$ and  for  $K>0$ the constant in p. 249
  \cite{WW};}\\ & \Gamma \left (\varsigma   \right )= \frac{\Gamma \left (\varsigma  +2 \right )}{\varsigma \left (\varsigma  +1 \right )} . \end{aligned}
\end{equation*}
Then the absolute value of the integral is bounded by
\begin{equation*}  \begin{aligned}  C_1  |z|^{-\frac 32}&  \int _{\R }  e ^{-t\arg (z) }
  e^{-t \arg \left ( \frac{1}{2} +\im t  \right ) } \left |\frac{3}{2}-\im  t \right | ^{-1} \left |\frac{1}{2}-\im  t \right | ^{-1}\times \\&
  \times e^{\left (t + \frac{\nu _0}{2}\right ) \arg \left (2-\im \left (t + \frac{\nu _0}{2}\right ) \right ) }
  \left |2-\im \left (t + \frac{\nu _0}{2}\right ) \right | ^{\frac 32}
  e^{\left (t + \frac{\nu _0}{2}\right ) \arg \left (\frac{3}{2}-\im \left (t + \frac{\nu _0}{2}\right ) \right ) }
 \left |\frac{3}{2}-\im \left (t + \frac{\nu _0}{2}\right ) \right |
   dt \\    \le C_2  |z|^{-\frac 32}&   \int _{\R }  e ^{-t\arg (z)-|t| \frac{3}{2}\pi  } \langle t \rangle ^{\frac{1}{2}}  dt\le C_3  |z|^{-\frac 32}  \left (\frac{3}{2}\pi -|\arg (z)|\right )^{- \frac{3}{2}}   \end{aligned}
\end{equation*}
for fixed constants  which depend on $\rho _0$
and for $|\arg (z)| <\frac{3}{2}\pi $.  This  and the identity \eqref{eq:idD}   yield   inequality  \eqref{eq:p21}
if $\zeta$ is outside a union of preassigned small cones containing $\Sigma _P$. Near the cones we can proceed  by estimating similarly
the r.h.s.'s of the identities
\begin{equation*}  \begin{aligned} &    D_{\im \nu _0}(\zeta) =e^{-\nu _0 \pi  }   D_{\im \nu _0}(-\zeta) +
 \frac{\sqrt{2\pi}}{\Gamma (-\im \nu _0)}   e^{\frac  \im 2 (\im \nu _0 +1)   \pi  }     D_{-\im \nu _0-1}(-\im \zeta)   \ , \\&      D_{\im \nu _0}(\zeta) =e ^{\nu _0 \pi  }   D_{\im \nu _0}(-\zeta) +
 \frac{\sqrt{2\pi}}{\Gamma (-\im \nu _0)}    e^{-\frac  \im 2 (\im \nu _0 +1)   \pi  } D_{-\im \nu _0-1}( \im \zeta) \, . \end{aligned}
\end{equation*}
This completes the proof of Lemma \ref{lem:asP}. \qed

\subsubsection{The $\overline{\partial}$ argument}
\label{ssec:dbar}
We follow closely the argument of
Dieng and McLaughlin \cite{DM} which have a simpler discussion than in \cite{DZ2,DZ,DP}   as to  how to localize the RH to the model RH problem.  We modify slightly   \cite{DM}
to allow the case $s\in (1/2, 1)$ in Theorem \ref{th:purerad}.

We fix a smooth cutoff function of compact support, with $\chi (x)\ge 0$   for any $x$ and $\int \chi dx =1$. For $ \varepsilon
\neq 0$ let $ \chi  _{\varepsilon}(x)=  \varepsilon ^{-1} \chi (\varepsilon ^{-1}x).$
   For $z\in \C$ and for the convolution $f*g(x)=\int f(x-y) g(y) dy$, we define $\mathbf{r}(z)$ as follows:

   \begin{equation} \label{eq:defr} \mathbf{ r}(z) =
\left\{\begin{matrix}
          r(\Re z) \text{ for $\Im z =0$} \,\\
        \chi  _{\Im z} * r  (\Re z)  \text{  for $\Im z \neq 0$,}
\end{matrix}\right.
\end{equation}

The first  step is    the following proposition.

\begin{proposition}\label{prop:DM1}  Set $\widehat{r} _0= r(z_0) e^{-2\im \nu (z_0) -2\beta (z_0, z_0)}$ as in \eqref{eq:defVP1}.
Fix $\lambda _0>0$ and assume $\| r\| _{H ^{s }}<\lambda _0$  for a preassigned $s\in (1/2, 1]$.
Then there exist      functions $R_j$ defined  in $\overline{\Omega} _j$ for $j=1,3,4,6$ and a constant  $c   $
such that the following properties hold:
\begin{equation*}  \begin{aligned}   \left\{\begin{matrix} R_1(z)=r(z ) \text{ for  }   z-z_0\in \R _+ , \\  R_1(z)=f_1(z-z_0):=\widehat{r} _0  (z-z_0) ^{-\im \nu (z_0)}   \delta ^{2}(z )\text{ for  }   z-z_0\in e^{\im  \frac{\pi}{4}}   \R _+  ; \end{matrix}\right.
\end{aligned}
\end{equation*}

\begin{equation*}  \begin{aligned}   \left\{\begin{matrix} R_3(z)=\frac{\overline{r}(z )}{1+|{r}(z )|^2} \text{ for  }   z-z_0\in \R _- , \\  R_3(z)=f_3(z-z_0):=  \frac{\overline{\widehat{r} _0} }{1+|{r}(z_0 )|^2}    (z-z_0) ^{ \im \nu (z_0)}   \delta ^{-2}(z )\text{ for  }   z-z_0\in e^{3\im \frac{\pi}{4}}   \R _+  ; \end{matrix}\right.
\end{aligned}
\end{equation*}

\begin{equation*}  \begin{aligned}   \left\{\begin{matrix} R_4(z)=\frac{ {r}(z )}{1+|{r}(z )|^2} \text{ for  }   z-z_0\in \R _- , \\  R_4(z)=f_4(z-z_0):=  \frac{ {\widehat{r} _0} }{1+|{r}(z_0 )|^2}    (z-z_0) ^{- \im \nu (z_0)}   \delta ^{ 2}(z )\text{ for  }   z-z_0\in e^{5\im \frac{\pi}{4}}   \R _+  ; \end{matrix}\right.
\end{aligned}
\end{equation*}

\begin{equation*}  \begin{aligned}   \left\{\begin{matrix} R_6(z)=\overline{r}(z ) \text{ for  }   z-z_0\in \R _+ , \\  R_6(z)=f_6(z-z_0):=\overline{\widehat{r}} _0  (z-z_0) ^{- \im \nu (z_0)}   \delta ^{-2}(z )\text{ for  }   z-z_0\in e^{-\im \frac{\pi}{4}}   \R _+ ; \end{matrix}\right.
\end{aligned}
\end{equation*}

$\forall$ $j\in \{1,3,4, 6\} $, $\forall$   $z\in  \Omega _j+z_0$  and for $ \varphi (x) = - \chi (x) - x\chi ' (x) $, we have for a fixed $c$
\begin{equation} \label{eq:estdbar} \begin{aligned}      |\overline{\partial}R_j(z)|\le c \| r\| _{H^s(\R)} |z-z_0|^{s-\frac{3}{2}}+ c  | \partial _{\Re z}\mathbf{ r}(z)| +c |  (\Im z) ^{-1}  \varphi  _{\Im z} * r    (\Re z)|
\end{aligned}
\end{equation}

\end{proposition}
\proof
  The  $R_j(z)$  can be defined explicitly. For $j=1,3$ in particular, we  set
  for $z-z_0=u+i\nu $ and $b(x )=  \cos (2x)$,

\begin{equation}\label{eq:r1}  \begin{aligned}   R_1(z)&= b (  \arg (u+i\nu) )\mathbf{r} (  z) + (1-b (  \arg (u+i\nu) ))f_1(u+i\nu ),\\   R_3(z)&= \cos (2 (\arg (z-z_0)  -\pi ) )\frac{\overline{\mathbf{r}}(  z )}{1+|{\mathbf{r}}(  z )|^2} \\& + (1-\cos (2 (\arg (z-z_0)  -\pi ) ))f_3(u+i\nu ) .
\end{aligned}
\end{equation}
The other $R_j(z)$'s can be defined similarly.  This yields functions with the desired boundary values.  Now we prove the bounds, and for definiteness we consider case $j=1$ only.  We have
$$\overline{\partial} R_1  = (\mathbf{r} -f_1) \overline{\partial}  b + \frac{b}{2}  (  \chi  _{\Im z} * r '  (\Re z)  +\im  (\Im z) ^{-1} \varphi  _{\Im z} * r    (\Re z) ), $$
with $\varphi (x) = - \chi (x) - x\chi ' (x) $. Notice that $\widehat{\varphi} (0)=0$.
Then we have the bound

\begin{equation*}   \begin{aligned}   |\overline{\partial} R_1| &\le |\chi  _{\Im z} * r '  (\Re z) |+| (\Im z)  ^{-1} \varphi  _{\Im z} * r    (\Re z)| \\&+ \frac{c}{|z-z_0|} \left ( | r(z)-r(z_0)| + | f_1(z)-r(z_0)|   \right ) .
\end{aligned}
\end{equation*}
To obtain the desired estimate for $ |\overline{\partial} R_1|$ we need to bound the last line. By \eqref{eq:beta3} we have $| r(z)-r(z_0)| \le C |z-z_0| ^{s-\frac{1}{2}}\|  r \| _{H^s}.$
Next, we have
\begin{equation*}   \begin{aligned} & f_1(z)-r(z_0) = r(z_0)  \times \\& \left [ \exp \left( 2\im \nu (z_0) ( (z-z_0) \log (z-z_0) -(z-z_0+1) \log (z-z_0+1) )+ 2(\beta (z, z_0) - \beta (z_0, z_0)  ) \right )-1  \right ]    .
\end{aligned}
\end{equation*}
By Lemma \ref{lem:beta} we have $\beta (z, z_0) - \beta (z_0, z_0)= C(\rho ,s)\| r\| _{H^s}|z-z_0| ^{s-\frac{1}{2}}$. Since for $z$ close  to $z_0$ both $(z-z_0) \log (z-z_0)$  and $(z-z_0+1) \log (z-z_0+1)$ are   $O(|z-z_0| ^{s-\frac{1}{2}})$, we get the  desired estimate for $ |\overline{\partial} R_1|$.

 \qed

We now extend as follows the  matrices in \eqref{eq:jump1}:
\begin{equation} \label{eq:jump2} \begin{aligned}  &
 W_R =  \begin{pmatrix}
    1   & 0 \\  e^{ 2\im t \theta } R_1
  & 1
 \end{pmatrix}     \text{ in } \Omega_1+z_0
  , \
 U_R = \begin{pmatrix}
    1   & e^{-2\im t \theta } R_3   \\   0
  & 1
 \end{pmatrix}
  \text{ in } \Omega_
 3+z_0  , \\&
  U_L=  \begin{pmatrix}
    1   & 0     \\   e^{ 2\im t \theta }  R_4
  & 1
 \end{pmatrix}  \text{ in } \Omega_4 +z_0    , \,
W_L=
\begin{pmatrix}
    1   & e^{ -2\im t \theta }   R_6   \\  0
  & 1
 \end{pmatrix}     \text{ in } \Omega_6+z_0      . \end{aligned}
\end{equation}
We set
\begin{equation}\label{eq:eqA}A:=
\left\{  \begin{matrix}
   m W^{-1}_R       \text{  in $\Omega _1+z_0$}   ,
 \\   m      \text{  in $(\Omega _2\cup \Omega _5)+z_0$}   ,\\
 m U^{-1}_R \text{  in $\Omega _3+z_0$},
 \\   m U_L   \text{  in $\Omega _4+z_0$},   \\ m W_L \text{ in $\Omega _6+z_0$}.
     \end{matrix}\right.
\end{equation}
 We set  $B :=A\delta ^{-\sigma _3}$,  obtaining a new  function with jump relations
   $B_+(z)=B_-(z)V_B (z)$     with jump matrix defined by
\begin{equation*}V_B (z):=
\left\{  \begin{matrix}
   \begin{pmatrix}
    1   & 0 \\  e^{ 2\im t \theta } R _1(z) \delta ^{-2}(z)
  & 1
 \end{pmatrix}      \text{ for }    z \in  z_0+e^{\im  \pi /4 }\R _+,
 \\  \begin{pmatrix}
    1   & - e^{-2\im t \theta } R_3 (z)\delta ^{-2}(z) \\  0
  & 1
 \end{pmatrix}  \text{ for }    z \in  z_0+e^{3 \im  \pi /4 }\R _+,\\
 \begin{pmatrix}
    1   & 0   \\   e^{ 2\im t \theta }  R_4 (z)\delta ^{-2}(z)
  & 1
 \end{pmatrix}\text{ for }    z \in  z_0+e^{5\im \pi /4 }\R _+,
 \\  \begin{pmatrix}
    1   & -
    e^{-2\im t \theta } R_6 (z)\delta ^{ 2}(z) \\  0
  & 1
 \end{pmatrix}  \text{ for }    z \in  z_0+e^{  -\im \pi /4  }\R _+.
     \end{matrix}\right.
\end{equation*}
Set  now $E(z):=B(z)P^{-1}(\sqrt{8t}(z-z_0))$.   By the choice \eqref{eq:defVP1} of the parameter $r_0$ in \eqref{eq:defVP}, 
    the jump
 matrices of $B(z)$ and of $P (\sqrt{8t}(z-z_0))$ coincide. This is elementary to check and holds for the same reasons of \cite{DM}.
 As a consequence, $E(z)$   does not have jump discontinuities.
 We now reverse the construction, we define $E$ using Corollary \ref{prop:dbar1} below
 and define  $B(z)$  by
 $B(z)=E (z) P (\sqrt{8t}(z-z_0))$. First though, we have the following auxiliary lemma, see \cite{DM}.

 \begin{lemma}\label{lem:lemJ} Let $\|  r\| _{H^s}\le \lambda _0$    for a preassigned $s\in (1/2, 1]$.
 Consider  the    following   operator
\begin{equation}  \label{eq:lemJ1} \begin{aligned} &    JH(z):= \frac{1}{\pi} \int _{\C}\frac{H( \varsigma )W( \varsigma )}{\varsigma -z}dA(\varsigma )
\end{aligned}
\end{equation}
with,      for $\zeta =\sqrt{8t}(z-z_0)$,\begin{equation*}W (z):=
\left\{   \begin{matrix}
   P(\zeta )\begin{pmatrix}
    0   & 0 \\  e^{ 2\im t \theta }\delta ^{-2}(z) \overline{\partial}R _1(z)
  & 0
 \end{pmatrix} P^{-1} (\zeta )    \text{ for }    z \in  \Omega _1,
 \\ P (\zeta ) \begin{pmatrix}
    0   & - e^{-2\im t \theta } \delta ^{-2}(z) \overline{\partial}R_3 (z)\\  0
  & 0
 \end{pmatrix} P^{-1} (\zeta )  \text{ for }    z \in  \Omega _3,\\ P (\zeta )
 \begin{pmatrix}
    0   & 0   \\   e^{ 2\im t \theta }  \delta ^{-2}(z)\overline{\partial}R_4 (z)
  & 0
 \end{pmatrix}P^{-1}(\zeta )  \text{ for }    z \in  \Omega _4,
 \\  P(\zeta )\begin{pmatrix}
    0   & -
    e^{-2\im t \theta } \delta ^{ 2}(z)\overline{\partial}R_6 (z) \\  0
  & 0
 \end{pmatrix} P^{-1}(\zeta )   \text{ for }    z \in  \Omega _6,\\ 0 \text{ for }    z \in  \Omega _2\cup \Omega _5.
     \end{matrix}\right.
\end{equation*}
Then,  we have  $J:L^ \infty (\C )\to L^ \infty  (\C )\cap C^0 (\C )$ and
  there exists a    $C=C (\lambda _0)$  s.t.
  \begin{equation}  \label{eq:lemJ2} \begin{aligned} &   \| J \| _{L^ \infty  (\C ) \to L^ \infty  (\C )}\le C t^{\frac{1-2s}{4}} \text{ for all $t>0$}.
\end{aligned}
\end{equation}
 \end{lemma}
\proof     For definiteness let $H\in L^\infty   (\Omega _1)$. Then
\begin{equation}  \label{eq:lemJ-3}    \begin{aligned} &    \pi |JH (z)| \le  \| H\| _{L^\infty }   \| \delta ^{-2}\| _{L^\infty (\Omega _1)} \int _{\Omega _1}
  \frac{ |\overline{\partial}R_1(\varsigma )  e^{2\im t \theta }|  }{| \varsigma  -z|}dA(\varsigma ).
\end{aligned}
\end{equation} We have $\| \delta ^{-2}\| _{L^\infty (\Omega _1)} \le 1$ by Prop. \ref{prop:scalarRH}.
By \eqref{eq:estdbar}  to   bound   \eqref{eq:lemJ-3} it is enough to bound $I_j$ for $j=1,2,3$ with
\begin{equation}  \label{eq:lemJ3}  \begin{aligned} &  I_j= \int _{\Omega _1}
  \frac{ |X_j (\zeta )  e^{2\im t \theta }|  }{|\varsigma -z|}dA(\varsigma ),  \quad X_1 (z ):= \partial _{\Re z}\mathbf{ r}(z)
   , \\&   X_2 (z ):= \| r\| _{H^s(\R)} |z-z_0|^{s-\frac{3}{2}} ,  \quad X_3 (z ):=    (\Im z) ^{-1} \varphi  _{\Im z} * r    (\Re z) .
\end{aligned}
\end{equation}
The estimates are like  those  in Sect. 2.4 \cite{DM}. We have, for $\varsigma-z_0=u+\im \nu $
and for $z-z_0=\alpha +\im \beta $,
\begin{equation}  \label{eq:lemJ33}    \begin{aligned} &  I_1=  \int _{\Omega _1}
  \frac{ | \partial _{u}\mathbf{ r}(\varsigma )| e^{-8t u\nu }   }{|\varsigma  -z|}du d\nu
  \le  \int _0 ^\infty   d\nu \int _{\nu}^{\infty}  \frac{ | \partial _{u}\mathbf{ r}(\varsigma )| e^{-8t u\nu }   }{|\varsigma -z|}du
  \\&
  \le \int _0 ^\infty   d\nu  e^{-8t  \nu ^2 }\|  \partial _{u}\mathbf{ r}( u , \nu )\| _{L^2_u}  \|   ( (u-\alpha )^2 +(\nu -\beta  )^2 ) ^{-1}\| _{L^2_u ( \nu  , \infty )} .
\end{aligned}
\end{equation}
By elementary computation we have $\|   ( (u-\alpha )^2 +(\nu -\beta  )^2 ) ^{-1}\| _{L^2_u ( \nu  , \infty )} \le C | \nu -\beta | ^{-\frac 12 }$, see \eqref{eq:I212} below. By Plancherel we have for fixed $C$
\begin{equation} \label{eq:I1}\begin{aligned} &
 \|  \partial _{u}\mathbf{ r}( u , \nu )\| _{L^2_u}  =  \|  \partial _{u} \int _{\R}\nu ^{-1}
 \chi ( \nu ^{-1}  (u-t))  r(t) dt \| _{L^2_u} = \|   \xi
 \widehat{\chi} ( \nu \xi )  \widehat{r} (   \xi )  \| _{L^2 }\\& \le  \nu ^{s-1}  \|  \xi ^{1-s}
 \widehat{\chi} (  \xi ) \| _{L^\infty } \| r \| _{H^s } \le C \nu ^{s-1}  \| r \| _{H^s } .
\end{aligned}
\end{equation}
 So
 \begin{align}\nonumber &  I_1\le C \| r \| _{H^s }   \int _\R   d\nu  e^{-8t  \nu ^2 } |\nu| ^{s-1}  | \nu -\beta | ^{-\frac 12 } \le  C t^{\frac{1-2s}{4}} \| r \| _{H^s }   \int _\R   d\nu  e^{-8   \nu ^2 } (  |\nu| ^{s-\frac 32} +   | \nu -\sqrt{t}\beta | ^{s-\frac 32} )\\& \le  (3C    \int _\R    e^{-8   \nu ^2 }    |\nu| ^{s-\frac 32} d\nu ) \| r \| _{H^s } t^{\frac{1-2s}{4}}
      .\label{eq:I11}
\end{align}
For the last inequality we used  the fact that for any $c\in \R$
 \begin{align} \nonumber&   \int _\R     e^{-8   \nu ^2 }   | \nu -c | ^{s-\frac 32}  d\nu  \le    \int _{|\nu |\le | \nu -c | }    e^{-8   \nu ^2 }   | \nu   | ^{s-\frac 32}  d\nu +     \int _{|\nu |\ge  | \nu -c | }    e^{-8   (\nu -c) ^2 }   | \nu -c | ^{s-\frac 32} d\nu\\& \le 2  \int _\R    e^{-8   \nu ^2 }    |\nu| ^{s-\frac 32} d\nu
      .\label{eq:tool}
\end{align}
The estimate for $I_3$ is similar  after replacing \eqref{eq:I1} with
\begin{equation} \label{eq:I3}  \begin{aligned} &
    \|   \nu ^{-2} \int
 \varphi  ( \nu ^{-1}  (u-t))  r(t) dt \| _{L^2_u} = \|   \nu ^{-1} \xi ^{-s}
 \widehat{\varphi} ( \nu \xi ) \xi ^{ s}  \widehat{r} (   \xi )  \| _{L^2 }\\& \le  \nu ^{s-1}  \|  \xi ^{ -s}
 \widehat{\varphi } (  \xi ) \| _{L^\infty } \| r \| _{H^s } \le C \nu ^{s-1}  \| r \| _{H^s },
\end{aligned}
\end{equation}
where the latter bound holds since $\widehat{\varphi }$ is a fixed Schwartz function with $\widehat{\varphi }(0)=0$.  Proceeding like in \eqref{eq:lemJ33}, we finally consider
\begin{equation} \label{eq:I2}   \begin{aligned} &  I_2\le    \int _0 ^\infty  e^{-8t  \nu ^2 }  d\nu     \|  |\varsigma  -z_0|^{s-\frac{3}{2}}  \| _{L^p(\nu ,\infty )}  \| |\varsigma -z |^{-1}  \| _{L^q(\nu ,\infty )}
\end{aligned}
\end{equation}
with an appropriate pair $ 1/p+1/q=1$. By \cite{DM} we have
\begin{equation} \label{eq:I212}   \begin{aligned} & \| |\varsigma -z |^{-1}  \| _{L^q(\nu ,\infty )}\le C |\nu  -\beta| ^{\frac{1}{q} -1}
\end{aligned}
\end{equation}
  and
\begin{equation} \label{eq:I22}  \begin{aligned} &
   \|  |\varsigma  -z_0|^{s-\frac{3}{2}}  \| _{L^p(\nu ,\infty )} = (\int _\nu ^\infty  |u+\im \nu | ^{ p(s-\frac{3}{2})} du ) ^{\frac{1}{p}}   = (\int _\nu ^\infty  (u ^2+ \nu  ^2 ) ^{ p \frac{2s- {3}}{4} }  du  ) ^{\frac{1}{p}} \\& =\nu  ^{   \frac{2s- {3}}{2} + \frac{1}{p} } (\int _1 ^\infty  (u ^2+1 ) ^{ p \frac{2s- {3}}{4} }  du  ) ^{\frac{1}{p}}.
\end{aligned}
\end{equation}
So by \eqref{eq:I2} and using again \eqref{eq:tool},  we obtain
\begin{equation} \label{eq:I2bis}   \begin{aligned} &  I_2\le C '  \int _0 ^\infty  e^{-8t  \nu ^2 }    \nu  ^{   \frac{2s- {3}}{2} + \frac{1}{p} }  |\nu  -\beta| ^{\frac{1}{q} -1}  d\nu \le 3 C '  \int _0 ^\infty  e^{-8t  \nu ^2 }    \nu  ^{   \frac{2s- {3}}{2}  }   d\nu \le C  t^{\frac{1-2s}{4}} .
\end{aligned}
\end{equation}
 The proof that $J(L^\infty )\subset C^0$ can be seen by the above estimates using
  standard facts, like dominated convergence, and is skipped here.

\qed

Taking  $E$   as solution of $E=1+J(E)$ we obtain the following result.

 \begin{corollary}\label{prop:dbar1}Fix $\lambda _0>0$ and assume $\| r\| _{H ^{s }}<\lambda _0$. Then there exist a  constant     $T $   such that for $t\ge T$ there exists a $E(z)$ continuous in $\C$  and  satisfying the following additional properties:

\begin{itemize}

\item  [(1)] $E(z)$ is continuous in $\C$,

\item  [(2)] $E$ solves the system  $\overline{\partial}E=EW$,

\item   [(3)] $E(z)\to 1$  for $z\to \infty$.

\end{itemize}

 \end{corollary}

 \qed

Claim (3) in {Corollary} \ref{prop:dbar1}  can be replaced by the following
sharper result.

 \begin{lemma}\label{lem:E1}  There exists $ \varepsilon _0>0$ such that for
   $\| r\| _{H ^{s }}<\varepsilon _0$  there exist  constants    $T $ and $c $  such that for $t\ge T$ and for $z\in \Omega _2 \cup \Omega _5$
\begin{equation} \label{eq:dbar11} \begin{aligned}&   E(z) =1+ \frac{E_1}{z} + O(   z ^{-2})\\& |E_1|\le c \| u_0\| _{L ^{2,s}} t^{-\frac{2s+1}{4}} \text{ for } t\ge T.
\end{aligned}
\end{equation}
 \end{lemma}
\proof   We have $E_1= \frac{1}{\pi }  \int _\C EW dA ,$ so  $|E_1|\le \frac{\| E\| _\infty}{\pi }   \sum _j\int _{\Omega _j}  |W| dA .$ We bound the integrals using a decomposition as in
  \eqref{eq:lemJ3} and for definiteness we consider only $j=1$. For $\ell =1,3$ we have
  by \eqref{eq:I1} and  \eqref{eq:I3} and starting as in \eqref{eq:lemJ33},
 using
 \begin{equation*}
    \|   e^{-8t u\nu }  \| _{L^q_u ( \nu  , \infty )}= ( \int _\nu ^\infty    e^{-8qt u\nu } du ) ^{\frac{1}{q}} =(8qt \nu ) ^{-\frac{1}{q}}e^{-8qt  \nu ^2 } ,
 \end{equation*}
 we have
  \begin{equation} \label{eq:dbar111} \begin{aligned}&  \int _{\Omega _1}  |X_\ell (\zeta )  e^{2\im t \theta }| dA \le \| r \| _{H^s } \int _{0}^{\infty} \nu ^{s-1}  \|   e^{-8t u\nu }  \| _{L^2_u ( \nu  , \infty )}d\nu \\& \le C' t ^{-\frac{1}{2}}\int _{0}^{\infty} \nu ^{s-\frac{3}{2}}  e^{- t  \nu ^2 }d\nu  \| r \| _{H^s }  = C_s  t^{-\frac{2s+1}{4}}  \| r \| _{H^s }.
\end{aligned}
\end{equation}
 For $\ell =2$  we  use \eqref{eq:I22}   and  the elementary bound
\begin{equation} \label{eq:dbar112} \begin{aligned}&  \int _{\Omega _1}  |X_2 (\varsigma  )  e^{2\im t \theta }| dA \le \| r \| _{H^s } \int _{0}^{\infty}   \|  | \varsigma-z_0|^{s-\frac{3}{2}}  \| _{L^p(\nu ,\infty )}   \|   e^{-8t u\nu }  \| _{L^q_u ( \nu  , \infty )}d\nu \\& \le C  \| r \| _{H^s }    t  ^{  -\frac{1}{q}} \int _{0}^{\infty} \nu  ^{   \frac{2s- {3}}{2} + \frac{1}{p} -\frac{1}{q }  }   e^{-t\nu^{2}} d\nu \le   C_s  t^{-\frac{2s+1}{4}}  \| r \| _{H^s }.
\end{aligned}
\end{equation}
 Then we get  \eqref{eq:dbar11} by $\| r \| _{H^s } \le C \| u_0\| _{L ^{2,s}}$
 for a fixed $C$ by the Lipschitz continuity of Lemma \ref{lem:dirsc2}.

\qed

  Theorem \ref{th:purerad} follows by        $m(t,x,z)=E(z) P(\sqrt{8t}(z-z_0)) \delta ^{\sigma_3}(z) $    in $\Omega _3 +z_0$,   with
\begin{equation*}      \begin{aligned} &   m (z) =1+\frac{m_1}{z}+O(z^{-2}) \text{ with }  m_1=E_1+\frac{P_1}{\sqrt{8t}}+\begin{pmatrix}
\delta _{1}   & 0 \\
0 & -\delta _{1}
 \end{pmatrix} ,
\end{aligned}
\end{equation*}
where by \eqref{eq:invscat1} and Proposition \ref{prop:dbar1}
for $t\ge T(s,\lambda _0)  $ and a fixed  $C=C(s,\lambda _0) $  we have
\begin{equation*}      \begin{aligned} &   |u(t,x)-   2\im \frac{k_1}{\sqrt{8t}} |\le C |t| ^{-\frac{2s+1}{4}} \text{ and }  |u(t,x)  |\le C |t|^{-\frac{1}{2}},
\end{aligned}
\end{equation*}
where we recall that we have fixed $s\in (1/2 , 1]$.
The time reversibility of  the
NLS \eqref{eq:nls}  (see also later in   {Lemma}
  \ref{lem:distgs}) yields  the  same estimates   also   $\forall$ $t\le - T(\lambda _0)$. This proves Theorem \ref{th:purerad}
for $u_0\in H^1 (\R )\cap  L^{2,s}(\R )$.

Consider     $u_0\in L^{2,s}(\R )$   but  $u_0\not \in H^1 (\R )$.
Let  $  u(t) $ be the solution,   provided by  Theorem \ref{th:L2}, of the corresponding Cauchy problem \eqref{eq:nls}.
Consider   a sequence $u_{0 n}\in L^{2,s}(\R) \cap H^1(\R)$
such that
$ u_{0n} {\rightarrow}u_{0 } $ in $L^{2,s}(\R)$.  Then   for the reflection coefficients  we have
$ r_{ n} {\rightarrow}r $ in $H^{s}$ by Lemma \ref{lem:dirsc2}.

We can assume $\|   r_{ n} \| _{H^{s}}
\le 2 \|    r\| _{H^{s}}  $ for all $n$.   By the discussion developed so far,
there is a fixed $C$, which depends  only on  $\lambda _0$, where $\lambda _0 \ge \|    r\| _{H^{s}} $,
such that  for $|t|\ge T(\lambda _0) $  we have $ |u_{n}(t,x)|\le C |t| ^{-\frac{1}{2}}$ for almost any $x$.
By  Theorem \ref{th:L2} we known that for any $t$ we have  $ u_{n}(t) {\rightarrow}u(t) $ in
$L^2 (\R )$. This implies that  for almost any $x$ we have $ u_{n}(t,x) {\rightarrow}u(t,x) $. In turn,
we can conclude that   $ |u (t,x)|\le C |t|  ^{-\frac{1}{2}}$ for almost any $x$. This completes the proof  of the statement in
Theorem \ref{th:purerad} also in the case when  $u_0\in L^{2,s}(\R )$   but  $u_0\not \in H^1 (\R )$.

\subsubsection{Several remarks }
\label{ssec:remarks}

  Lemma \ref{lem:lemJ} yields   $\| E-1\| _{L^\infty (\C )} \le C  t^{\frac{1-2s}{4}}$. However
we will need the following lemma.

\begin{lemma}
  \label{lem:dbar2}  Let   $z_1\in \C _+$.    Assume  $  \|    u_0 \| _{L^{2,s}(\R) } < \varepsilon _0$.   Then
 there   are a  $\varepsilon _0>0$, a
      $c>0$ and a  $T>0$    such that
 \begin{equation} \label{eq:dbar14} \begin{aligned} &
 |1-E(z_1)|\le  c  \  t^{-\frac{2s+1}{4}} \|    u_0 \| _{L^{2,s}(\R) }   \text{ for } t\ge T  .
\end{aligned}\end{equation}

\end{lemma}
\proof
The argument is like in Lemma \ref{lem:E1}. We have $|E-1|\le  \frac{\| E\| _\infty}{\pi }   \sum _j\int _{\Omega _j}  \frac{|W|}{|\zeta -z_1|} dA $. Once again, we estimate only the term with $j=1$.   Using the notation in
\eqref{eq:lemJ3} and proceeding like in \eqref{eq:lemJ33}, for $ z_1=\alpha _1+\im \beta _1$ we have for $\ell=1,3$
\begin{equation} \label{eq:e111} \begin{aligned}&  \int _{\Omega _1}    \frac{|X_\ell (\varsigma )  e^{2\im t \theta }|}{|\varsigma -z_1|} dA \le \| r \| _{H^s } [ A_1+A_2] \, , \, A_1:=
\int _{0}^{\frac{\beta _1}{2}} \nu ^{s-1}  \|   \frac{e^{-8t u\nu }}{|\varsigma -z_1|}  \| _{L^2_u ( \nu  , \infty )}d\nu
\\& A_2:= \int _ {\frac{\beta _1}{2}}^{\infty} \nu ^{s-1}  \|   \frac{e^{-8t u\nu }}{|\varsigma -z_1|}  \| _{L^2_u ( \nu  , \infty )}d\nu .
\end{aligned}
\end{equation}
As in \eqref{eq:dbar111} we have
\begin{equation*}   \begin{aligned}&    A_1=\int _{0}^{\frac{\beta _1}{2}} \nu ^{s-1}  \|   \frac{e^{-8t u\nu }}{\sqrt{ (u+  \frac{x}{4t}-\alpha _1)^2+ (\nu-\beta _1)^2}}  \| _{L^2_u ( \nu  , \infty )}d\nu  \\& \le  C'(\beta _1) \int _{0}^{\frac{\beta _1}{2}} \nu ^{s-1}  \|  {e^{-8t u\nu }}   \| _{L^2_u ( \nu  , \infty )}d\nu \le C (s,\beta _1)t^{-\frac{2s+1}{4}}.
\end{aligned}
\end{equation*}
 By \eqref{eq:I212}, using  $t\ge 1$ and $e^{-8t  \nu ^2}\le e^{- t  \gamma_1 ^2}e^{-4  \nu ^2}$  for $\nu \ge \frac{\beta _1}{2}$,
and using bounds similar to those for \eqref{eq:I11}, we have
\begin{equation} \label{eq:e112} \begin{aligned}&    A_2\le \int _ {\frac{\beta _1}{2}}^{\infty}  e^{-8t  \nu ^2 }\nu ^{s-1}  \|     {|\varsigma -z_1|}^{-1}  \| _{L^2_u ( \nu  , \infty )}d\nu \le  C   \int _ {\frac{\beta _1}{2}}^{\infty}e^{-8t  \nu ^2 } \nu ^{s-1}  |\nu  -\beta _1| ^{-\frac{1}{2} }d\nu
\\& \le  C e^{- t  \beta_1^2 } \int _ {0}^{\infty}e^{-4   \nu ^2 } \nu ^{s-1}  |\nu  -\beta _1| ^{-\frac{1}{2} }d\nu \le C' e^{- t  \beta _1^2 }.
\end{aligned}
\end{equation}
Turning to the case $\ell =2$, we similarly have
\begin{equation*}   \begin{aligned}&  \int _{\Omega _1}    \frac{|X_2 (\varsigma )  e^{2\im t \theta }|}{|\varsigma -z_1|} dA \le \| r \| _{H^s } [ B_1+B_2] \, , \, \\& B_1:=
\int _{0}^{\frac{\beta _1}{2}}    \int _\nu ^\infty      |\varsigma -z_0|^{s-\frac{3}{2}}   \frac{e^{-8t u\nu }}{|\varsigma -z_1|}   d\nu \, , \,
B_2:= \int _ {\frac{\beta _1}{2}}^{\infty}     \int _\nu ^\infty      |\varsigma -z_0|^{s-\frac{3}{2}}   \frac{e^{-8t u\nu }}{|\varsigma -z_1|}   d\nu .
\end{aligned}
\end{equation*}
Then $B_1\le C (\beta _1, s) t^{-\frac{2s+1}{4}}$  by \eqref{eq:dbar112} and by $|\varsigma -z_1|\ge \beta _1/2$. We have
by \eqref{eq:I2}--\eqref{eq:I2bis} and using  $t\ge 1$
\begin{equation}   \begin{aligned}&    B_2\le \int _ {\frac{\beta _1}{2}}^{\infty}  e^{-8t  \nu ^2 }        \|  | \varsigma -z_0|^{s-\frac{3}{2}}  \| _{L^p(\nu ,\infty )}  \| |\varsigma -z _1|^{-1}  \| _{L^q(\nu ,\infty )}d\nu
\\& \le C e^{- t  \beta_1^2 } \int _ {0}^{\infty}  e^{-4  \nu ^2 }       \nu  ^{   \frac{2s- {3}}{2} + \frac{1}{p} }  |\nu  -\beta  _1| ^{\frac{1}{q} -1}  d\nu  \le C _s  e^{- t  y_1^2 }.
\end{aligned}
\end{equation}
\qed

\begin{lemma}
  \label{lem:pert}  Fix $z_1=\alpha _1+\im \beta _1$  with $\beta _1>0$. There is  $\varepsilon _0$ sufficiently small such that for
   $\|    u_0 \| _{L^{2,s}(\R) } < \varepsilon _0$
        there is a constant $C  $  such that
	
	\begin{equation} \label{eq:pert1} \begin{aligned} &  |    1- W_R  (z_1) |     \le C  e^{   - t 8\beta _1^2 } \|    u_0 \| _{L^{2,s}(\R) }  \text{  if }z_1\in \Omega _1+z_0\\&
    |    1-  U _R  ^{-1}(z_1) |     \le C  e^{   - t 8\beta _1^2 } \|    u_0 \| _{L^{2,s}(\R) }  \text{  if } z_1\in \Omega _3+z_0.
\end{aligned}
\end{equation}
	 \end{lemma}
\proof   By \eqref{eq:r1} we have that   $\| R _j\| _{L^\infty (\Omega _j+z_0)}\le  C'    \|  r \| _{H^s(\R)}\le C \|    u_0 \| _{L^{2,s}(\R) }$ for $j=1,3$.  If $ z_1\in \Omega _1+z_0$  we have $  \alpha _1+\frac x{4t}\ge  \beta _1$ and so
 $|e^{-2\im t\theta }|\le e^{ -8t(\alpha _1+\frac x{4t}) \beta _1}\le e^{   - t 8\beta _1^2 } $.   If $ z_1\in \Omega _3+z_0$  we have  similarly $|e^{ 2\im t\theta }|\le  e^{   - t 8\beta _1^2 } $.
This yields    \eqref{eq:pert1}.\qed

\begin{lemma}
  \label{lem:delta}  Fix $z_1= \alpha _1+\im \beta _1$    with $\beta _1>0$.
  Fix $\rho _0>0$. Let  $\rho :=\|  r \| _{L^\infty(\R)}$
  and assume $\rho<\rho _0$. Then there exists a constant   $C$ independent from  $z_0$  such that
  	\begin{equation} \label{eq:delta11} \begin{aligned} &   \left | \delta (z_1) - \Delta (z_1)
			\right | \le C   \|  r\| _{L^2}^2 \\& \text{where }   \Delta (z_1):=
			 \exp \left ( \frac 1{2\pi \im } \int _{-\infty}^{  \alpha _1}    \frac{\log (1+|{r} (\varsigma )|^2)}{\varsigma -z_1} d\varsigma  \right ).
\end{aligned}
\end{equation}
   Fix $K>0$.
  Then for $|z_0-\alpha _1|\le K  /\sqrt{t}$
 there exists a constant   $C$    such that
		
			\begin{equation} \label{eq:delta1} \begin{aligned} &   \left | \delta (z_1) - \Delta (z_1)
			\right | \le \frac{ C}{\sqrt{t}  \beta _1}   \log (1+ \rho ^2 ) .
\end{aligned}
\end{equation}
\end{lemma}

\proof     By Proposition \ref{prop:scalarRH},  we have  for a fixed $c$
\begin{equation*}
\begin{aligned} &    \left |  \gamma (z_1)-  \frac 1{2\pi \im } \int _{-\infty}^{  \alpha _1}    \frac{\log (1+|{r} (\varsigma )|^2)}{\varsigma -z_1} d\varsigma \right |  =  \frac 1{2\pi   } \left |\int _{  \alpha _1}^{z_0}    \frac{\log (1+|{r} (\varsigma )|^2)}{\varsigma- \alpha _1- \im \beta _1} d\varsigma\right  |  \le \frac{c }{\beta _1}     \|  r \| ^2_{L^2}    .
\end{aligned}
\end{equation*}
This yields \eqref{eq:delta11} since the bound $|\delta (z)|\le (1+ \rho ^2 ) $ is independent from $z_0$.
Similarly  \eqref{eq:delta1} follows from
\begin{equation*}
\begin{aligned} &    \left |  \gamma (z_1)-  \frac 1{2\pi \im } \int _{-\infty}^{ \alpha _1}    \frac{\log (1+|{r} (\varsigma )|^2)}{\varsigma -z_1} d\varsigma  \right |  =  \frac 1{2\pi   } \left |\int _{  \alpha _1}^{z_0}    \frac{\log (1+|{r} (\varsigma )|^2)}{\varsigma - \alpha _1- \im \beta _1} d\varsigma \right  | \le \frac{|z_0-\alpha _1| }{\beta _1}   \log (1+ \rho ^2 )  .
\end{aligned}
\end{equation*}
These yield Lemma \ref{lem:delta}. \qed

\vspace{0.25cm}

We will use the inequalities in Sect. \ref{ssec:remarks} for the proof of Theorem \ref{th:asstab}.
Notice that similar inequalities are also in Lemmas 5.18--5.21 \cite{DP}.

\section{Proof of Theorem \ref{th:asstab}}
\label{sec:asstab}

Recall  by Remark \ref{lem:soliton} that solitons (\ref{eq:groundstate})
belong to $ \mathcal{G}_1 $ (see under Lemma \ref{lem:dirsc1}). Since  $ \mathcal{G}_1 $ is an open subset of $L^1(\R )$,
see Lemma \ref{lem:dirsc1},   if the value of $\varepsilon _0>0$
in the bound \eqref{eq:asstab0} is small enough, then the initial datum
$u_0$ belongs to $ \mathcal{G}_1$. Notice also that the positive constant $\varepsilon_0$
can be taken  independent of $(\gamma _0,x_0)$.  Indeed, when we replace $u_0(x)$
with $u_0(x-x_0)$, their scattering function $a(z)$ is the same, while
$e^{\im \gamma _0} u_0(x)$   describes a compact set in $L^ {2,s}(\R )$   as $\gamma _0$ varies in $\R$.

We consider now an initial datum  $u_0$ satisfying the bound \eqref{eq:asstab0}.
The scattering datum   associated with the initial datum $u_0$,  which  by Lemma \ref{lem:dirsc1}
and Remark \ref{lem:soliton} belongs  to the space  $\mathcal{S}(1,1)$ defined in \eqref{eq:spdata},
 is close to those of the soliton
$\varphi_{\omega_0,\gamma_0,v_0}(0,x-x_0)$  by Lemmas \ref{lem:dirsc1} and \ref{lem:dirsc2}.
By  Lemma \ref{lem:dirsc2}, we know that  $u_0 \in L^{2,s}(\R)$ implies $r\in   H^s(\R)$.
Furthermore, by the Lipschitz continuity of $u_0\to r$ and the fact that the soliton
has $r \equiv 0$,  we have $\| r \| _{H^s(\R)}\le C \epsilon $,
with $C=C(\omega_0, v_0)$ and the value of $\epsilon $ is given in \eqref{eq:asstab0}.

We define now a map
\begin{equation}     \label{eq:darboux1}  \begin{aligned}
\mathcal{G}_1 \times \C _+ \times \C _* \ni (u_0,z_1,c_1) \mapsto \widetilde{u}_0 \in \mathcal{G}_0
\end{aligned}
\end{equation}
by means of the transformation
\begin{equation}  \label{eq:refcoeff}
\widetilde{r}(z) := r(z) \frac{z-z_1}{z-\overline{z}_1}.
\end{equation}
By its  definition,   $\widetilde{r} \in H^s(\R)$  if $r \in H^s(\R)$ and there is $C > 0$
such that $\|\tilde{r} \|_{H^s(\R)} \leq C \| r \|_{H^s(\R)}$.
We then define $\widetilde{u}_0 \in \mathcal{G}_0 \cap L^{2,s}(\R)$
by the reconstruction formula \eqref{eq:invscat1},
after the corresponding RH problem (i)--(iii) is solved for the scattering datum  in
$\mathcal{S}(1,0) = \{ \widetilde{r} \in H^s(\R)\}$, see \eqref{eq:spdata}. By Lemma \ref{lem:uc},
we know that   $\widetilde{u}_0 \in \mathcal{G}_0 \cap L^{2,s}(\R)$
with norm $\|  \widetilde{u}_0 \| _{L^{2,s}(\R)}\le C \| \widetilde{r}  \| _{H^s(\R)} \leq
C \epsilon$.

We now assume also that $u_0\in H^1(\R)$, to define the time evolution of the scattering data
in $\mathcal{S}(1,1)$ and $\mathcal{S}(1,0)$. Let
$\widetilde{{u}} (t )\in \mathcal{G}_0 \cap L^{2,s}(\R) \cap H^1(\R)$ be the solution of the
cubic NLS equation with the initial datum $\widetilde{u}_0$ and
$ {{u}} (t )\in \mathcal{G}_1  \cap L^{2,s}(\R) \cap H^1(\R)$ be the solution of the cubic
NLS equation with the initial datum $u_0$.

Denote the  solution of the RH problem (i)--(iv) associated to $\widetilde{u}(t)$
by $m(t,x,z)$. The two solutions $u(t)$  and   $\widetilde{u}(t)$ are related by
the auto--B\"{a}cklund transformation formula,  which we state now.

\begin{lemma} \label{lem:darboux}
We have
\begin{equation}  \label{eq:darboux2}
 {u} (t,x) =  \widetilde{{u}} (t,x) + \textbf{B}\, , \quad
\textbf{B} := 4 \Im (z_1)  \frac{\mathfrak{b}_1 \bar{\mathfrak{b}}_2}
{|\mathfrak{b}_1 |^2+| {\mathfrak{b}_2}|^2},
\end{equation}
where
\begin{eqnarray}
\label{b-1}
\mathfrak{b}_1 & := & e^{- \im x z_1}  m_{11} (t,x,z_1)   -  \frac{ c_1 m_{12} (t,x,z_1) e^{ \im xz_1 +4\im t {z}_1^2}   }{ 2\im \Im  (z_1)  },\\
\label{b-2}
\mathfrak{b}_2 & := & e^{- \im x z_1}  m_{21} (t,x,z_1)  - \frac{ c_1 m_{22} (t,x,z_1) e^{ \im xz_1 +4\im t {z}_1^2} }{  2 \im \Im  (z_1)}.
\end{eqnarray}
\end{lemma}

\proof
Note that $(\mathfrak{b}_1,\mathfrak{b}_2)^T$ is a solution of the spectral system \eqref{eq:zs},
and hence the transformation formula (\ref{eq:darboux2}) is a particular example
of the general auto--B\"{a}cklund transformation formula used in \cite{MP}
(after the transformation $\tilde{u} \to -\tilde{u}$ and $\mathfrak{b}_2 \to -\mathfrak{b}_2$,
which leaves (\ref{eq:zs}) invariant).
The particular expressions (\ref{b-1})--(\ref{b-2}) were used in \cite{DP}
and we give a sketch of the proof of this transformation formula from Appendix A in \cite{DP}.

We denote by $ {m}$ (resp. $\textbf{m}$)  the  solution of the RH problem (i)--(iv) associated to    $\widetilde{u}$   (resp. $u$). We set
 $ {\psi} = {m} e^{-\im \sigma _3 xz}$.
Then    consider the function $ \widehat{\psi} (x,z)$
\begin{equation*}  \begin{aligned} &   \widehat{\psi} (x,z):   =  \mathfrak{a} (x)
\mu (z)    \mathfrak{a}^{-1}(x)  {\psi} (x,z)  \mu  ^{-1}(z),
 \end{aligned}
\end{equation*}
where
$$
\mu (z):=
 \begin{pmatrix}
    z-z_1   & 0 \\  0
  &   z-\overline{z}_1
 \end{pmatrix}
$$
  and  $\mathfrak{a} =[\mathfrak{a}_1, \mathfrak{a}_2]$  with
\begin{equation*}  \begin{aligned} &
\mathfrak{a}_1(x):=  {\psi} (x,z_1)  \begin{pmatrix}
   1 \\       \frac {- {c}_1  }{ z_1-\overline{z}_1 }
 \end{pmatrix}    \, , \quad   \mathfrak{a}_2(x):=  {\psi} (x,z_1)  \begin{pmatrix}
     \frac {-\overline{c}_1  }{ z_1-\overline{z}_1 } \\    1
 \end{pmatrix} .
 \end{aligned}
\end{equation*}
By symmetries of the spectral system \eqref{eq:zs}  we have $\mathfrak{a}_2= \begin{pmatrix}
   0& -1 \\       1 & 0
 \end{pmatrix}     \overline{\mathfrak{a}}_1$.  Notice that
 $ \mathfrak{a}_1 =(\mathfrak{b}_1,\mathfrak{b}_2)^T$ is given by
 \eqref{b-1} and \eqref{b-2}.

The function  $ \widehat{\psi} (x,z)$ has poles only at  $z_1$ and $\overline{z}_1$;
     $ \widehat{m} (x,z) :=  \widehat{\psi} (x,z)e^{ \im \sigma _3 xz}$
satisfies \eqref{eq:resm}--\eqref{explicit-V-k}  for $k=1$.   Furthermore, $ \widehat{m} (x,z)$
satisfies (i)--(iv) of the RH  problem involving $ \widehat{V} _x(z)= e^{- \im \sigma _3 xz}(z) \widehat{V}(z)e^{  \im \sigma _3 xz} $  with
 \begin{equation*}
\widehat{V}(z) =   \widehat{\psi} _-^{-1}(x,z)\widehat{\psi}_+ (x,z)
=  \mu (z) {\psi} _-^{-1}(x,z) {\psi}_+ (x,z)  \mu ^{-1}(z)
= \begin{pmatrix}
    1+|\widehat{r}(z)| ^2  &   \overline{\widehat{r}  (z)} \\   {\widehat{r}} (z)  & 1
 \end{pmatrix},
\end{equation*}
where
$$
\widehat{r}(z) := r(z) \frac{z-\overline{z}_1  }{z- {z}_1}.
$$

All these formulas are in \cite{DP}, with a  different notation
(our  reflection coefficient $r(z)$ is equivalent to $\overline{r}(z)$ in
\cite{DP}, whereas our $z$ is $-z/2$ in  \cite{DP}).  It is clear by the uniqueness
of the inverse problem that $  {\textbf{m}}= \widehat{m}$.

We have expansions ${\textbf{m}} (x,z)=1 +\frac{{\textbf{m}}_1 (x )}{z}+o(z^{-1 })$
and  ${ {m}} (x,z)=1 +\frac{{ {m}}_1 (x )}{z}+o(z^{-1 })$.    By  an elementary computation,
we have $ \textbf{m}_1= {m}_1- \mathfrak{a}\mu _1\mathfrak{a} +\mu _1$, where
$$
\mu_1 := \begin{pmatrix}      z_1 & 0 \\  0 &\overline{z}_1 \end{pmatrix}.
$$
Therefore, the reconstruction formula (\ref{eq:invscat1}) yields
$$
u= \im [ \sigma _3,  \textbf{m}_1 ] _{12} =   \im [ \sigma _3,  {m}_1- \mathfrak{a}\mu _1\mathfrak{a} ] _{12},
$$
which  proves  \eqref{eq:darboux2}.\qed

\vspace{0.25cm}

\begin{remark}
The soliton in Remark \ref{lem:soliton} is obtained for $\widetilde{u} = 0$
and $m(x,z) = I$, when
$$
\mathfrak{b}_1 = e^{- \im x z_1} \quad \mbox{\rm and} \quad
\mathfrak{b}_2 = - \frac{c_1}{  2 \im \Im  (z_1)} e^{ \im xz_1 +4\im t {z}_1^2}.
$$
\end{remark}

By Theorem \ref{th:purerad}, we know that there exist constants $C_0 > 0$ and $T  > 0$
such that for all $|t| \geq T$, we have
$$
\|    \widetilde{u}(t, \cdot )  \| _{L^\infty ( \R)} \le    C_0 \epsilon |t|^{-\frac 12},
$$
since there is a constant $C > 0$ such that $\|  \widetilde{u}_0 \| _{L^{2,s}(\R)}\le C \epsilon$.

To prove  Theorem \ref{th:asstab} we need to focus only on  $\textbf{B}$.
>From the proof, we will see that the $(\omega _1 , v_1 )$  of the statement
of  Theorem \ref{th:asstab}  are those of the soliton with spectral data $(z_1,c_1)$.

We will consider only positive times, focusing on   $t\gg 1$. We know that
\begin{equation} \label{eq:mz1}m (t,x,z_1)=    \left\{\begin{matrix}
     E (z _1)  P(\sqrt{8t}  (z_1-z_0))  \delta ^{\sigma _3}  (z_1)  W_R(z_1)   \text{   if  }  z_1\in \Omega _1+ z_0 \ ,   \\
 E (z _1)  P(\sqrt{8t}  (z_1-z_0))  \delta ^{\sigma _3}  (z_1)      \text{   if  }  z_1\in \Omega _2+ z_0\ , \\
E (z _1)  P(\sqrt{8t}  (z_1-z_0))  \delta ^{\sigma _3}  (z_1)  U_R^{-1}(z_1)   \text{   if  }  z_1\in \Omega _3+ z_0 \  .
\end{matrix}\right.
\end{equation}
We have the following estimate.
\begin{lemma}\label{lem:est m}  Fix $\lambda _0>0$.  Then there is a $C>0$
and a $T>0$ such that for    $   \|  \widetilde{r} \|  _{H^s(\R)}< \lambda _0   $  we have for $t\ge T$
\begin{align} & \nonumber  |m    _{11}(t,x,z_1)  - \delta    (z_1)    |  +|m    _{22}(t,x,z_1)  - \delta  ^{-1}   (z_1)    |
  \le C    \|   \widetilde{r} \|  _{H^1(\R)} t^{-\frac 12}  (\|  \widetilde{r} \|  _{H^s(\R)}  +t^{- \frac{2s-1}4} )
\\&   \left |m    _{12}(t,x,z_1)  - \frac{\delta  ^{-1}   (z_1)   k_1}{\sqrt{8t} (z_1-z_0)}  \right |     +   \left |m    _{21}(t,x,z_1)  - \frac{\delta    (z_1)   k_2 }{\sqrt{8t} (z_1-z_0)}     \right |     \le C \|  \widetilde{r} \|  _{H^s(\R)}  t^{- \frac{2s+1}{4}}   .
  \nonumber \end{align}
\end{lemma}

\proof   By    Lemma  \ref{lem:dbar2}, we have $E (z _1)=1+O(\|  \widetilde{r} \|  _{H^s(\R)} t^{- \frac{2s+1}{4}})$.
 By    Lemma   \ref{lem:pert},  we have similar expansions for  $W_R(z_1) $  and $U_R(z_1)$.
We furthermore know by Proposition \ref{prop:scalarRH}   that
$|\delta  ^{\pm}   (z_1)|\le 1+\rho ^2$ for $  \rho =    \|  \widetilde{r} \| _{L^\infty(\R)}$.
>From Section \ref{ssec:model}, we recall the expansion
$$
P(\sqrt{8t} (z_1-z_0))= 1+ \frac{P_1}{\sqrt{8t} (z_1-z_0)}+O(\| \widetilde{r} \|  _{H^s(\R)} t^{-1}),
$$
where the $O$-term depends on  a fixed $ C=C(\lambda _0)$ and $P_1$  is given in  \eqref{eq:p1}.
We also recall that $|k_1|+|k_2|< C \|  \widetilde{r} \|  _{H^s(\R)}$.   These observations yield
Lemma \ref{lem:est m}. \qed

\vspace{0.25cm}

Now we start to analyze the term  $\textbf{B}$ in \eqref{eq:darboux2}.
Consider the following inequalities:
\begin{align} &  \label{eq:dom1}
|e^{- \im xz_1}  m _{11}(t,x,z_1) | > 10 \left | \frac{ {c}_1 m_{12} (t,x,z_1)
e^{ \im xz_1 +4\im t  {z}_1^2}   }{ 2\im \Im  (z_1)  }  \right |  \ , \\&
10 |e^{  -\im x z_1}   m_{21}(t,x,z_1) | <   \left | \frac{ c_1  m_{22} (t,x,z_1)
e^{ \im x z_1 + 4\im t z_1^2} }{  2 \im \Im  (z_1)}  \right |.  \label{eq:dom5}
\end{align}

\begin{lemma}
  \label{lem:bfB1}
Given    $\varepsilon  _0>0$ small,  there exist   $T(\varepsilon  _0)>0 $		 and $C>0$	
such that, if  $\| \widetilde{r} \| _{H^s(\R)}< \varepsilon  _0$ and	
if $(t,x)$ is such that at  least  one of   \eqref{eq:dom1}--\eqref{eq:dom5}  is false,
then we have $|\textbf{B}|<C t^{-\frac{1}{2}} \epsilon$ for $t\ge T(\varepsilon  _0)$.
\end{lemma}

\proof
Let us start by assuming that for  $(t,x)$ inequality
	\eqref{eq:dom1} is false. We are only interested to the
case when $t$ is large.  For the $\epsilon $  of \eqref{eq:asstab0}
and $\rho =\| \widetilde{r} \| _{L^\infty(\R)}$,
Lemma \ref{lem:est m}  implies for $t\ge T $,
\begin{eqnarray}
\nonumber |m_{12} | & \le &     (1+ \rho ^2)   |k_1|  t^{-\frac 12}     +   C \epsilon t^{-\frac {2s+1}4}    \\
\nonumber & \le &   t^{-\frac 12} \epsilon K \left(
 \frac{1}{2} (1+ \rho ^2) ^{-1}-     C  t^{ \frac { 1-2s}4} -C\epsilon \right ) \\
 & \le & t^{-\frac 12} \epsilon K  |m    _{22} |,
 \label{eq:dom111}
\end{eqnarray}
for a fixed and sufficiently large constant $K $.
Then,  if \eqref{eq:dom1} is false and  $t\ge T $,
both terms  in \eqref{eq:dom1}  are bounded from above by
\begin{equation*}
\left | \frac{ c_1  m_{22}(t,x,z_1) e^{ \im x z_1 + 4\im t z_1^2} }{  2 \im \Im  (z_1)} \right|.
\end{equation*}
For $t\ge T $ by the same argument of  \eqref{eq:dom111}
we have also
\begin{equation} \label{eq:dom2}
|e^{ -\im x z_1}  m _{21}(t,x,z_1) |  \le  t^{-\frac 12} \epsilon K |e^{- \im xz_1}  m _{11}(t,x,z_1) |.
\end{equation}
We conclude that  for $t\ge T $   and if  $(t,x)$ is in    the domain where \eqref{eq:dom1} is false,
we have for some fixed $K$
\begin{equation} \label{eq:dom3}
|\textbf{B}|     \le K
 \frac{\left |   m_{12}   e^{ \im xz_1 -4\im t  {z}_1^2}     \overline{m}_{22}
 e^{  -\im x\overline{z}_1 -4\im t  \overline{{z}}_1^2}  \right|}{|m_{22} e^{  \im x z_1 + 4\im t z_1^2}|^2}
 = K \frac{\left |   m_{12} \right |}{|m_{22}  | } \le \frac{C K}{\sqrt{t}}\epsilon.
\end{equation}
So now we assume that $(t,x)$  is such that    \eqref{eq:dom1} is true.
Notice that by  \eqref{eq:dom2}  and     \eqref{eq:dom1} we have for a fixed $K$
\begin{equation} \label{eq:dom4}
\frac{|\mathfrak{b}_1 e^{  \im x\overline{z}_1}   \overline{m_{21}}  |}
{\|\mathfrak{b}  \|^2 }  \le K
 \frac{|e^{- \im xz_1}  m _{11} e^{  \im x\overline{z}_1}   \overline{m_{21}}  |}
{|e^{- \im xz_1}  m _{11} |^2 } = K \frac{\left |   m_{21} \right |}{|m_{11}  | }
\le \frac{C K}{\sqrt{t}}\epsilon.
\end{equation}
Since we are  assuming that $(t,x)$ is such that   \eqref{eq:dom1}-- \eqref{eq:dom5}
are not both true, we assume now that
  \eqref{eq:dom1} is true and \eqref{eq:dom5}  is false.
   Then  by \eqref{eq:dom4}, for a fixed $K$
  \begin{equation} \label{eq:dom6}\begin{aligned} &  |\textbf{B}|\le 4|\Im z_1| \frac{|\mathfrak{b}_1 {\mathfrak{b}_2}|}
{\|\mathfrak{b}  \|^2 } \le K \frac{|\mathfrak{b}_1 e^{  \im x\overline{z}_1}   \overline{m_{21}}  |}{\|\mathfrak{b}  \|^2}
\le \frac{C K}{\sqrt{t}}\epsilon  .
   \end{aligned}
\end{equation}
The above inequalities prove  Lemma \ref{lem:bfB1} for values of $(t,x)$,
for which \eqref{eq:dom1}--\eqref{eq:dom5}  are not both true.\qed

\vspace{0.25cm}

We   assume  now that \eqref{eq:dom1}--\eqref{eq:dom5} are true. Then, by the last inequality in
\eqref{eq:dom3} and  by \eqref{eq:dom4}, up to terms bounded by $C  t^{-\frac{1}{2}}\epsilon  $,
what is left is the analysis of
\begin{equation}  \label{eq:dom7}
-2\im \frac{e^{- \im xz_1}  m _{11} \overline{c}_1  \overline{ m_{22}}
e^{ -\im x\overline{z}_1 -4\im t  \overline{{z}}_1^2}}{\|\mathfrak{b}  \|^2}.
\end{equation}
Set now
\begin{equation}   \begin{aligned} &
b^2 :=  |e^{- \im xz_1}  m _{11}  |^2+
\left | \frac{c_1  m_{22}   e^{ \im x z_1 + 4\im t z_1^2} }{  2 \im \Im  (z_1)}    \right |^2
     \nonumber
\end{aligned}\end{equation}
and expand
\begin{equation*}  \label{eq:dom8}  \begin{aligned} &
 \|\mathfrak{b}  \|^2 = b^2 \left  (1
  +  O\left (  b^{-1}  {\left | c_1  m_{12}   e^{ \im xz_1 +4\im t {z}_1^2}     \right |}   \right )  +
  O\left (  b^{-1}  {\left |  m_{21} e^{  -\im x z_1  }     \right |}   \right )    \right ).
\end{aligned}\end{equation*}
Then  the quantity in \eqref{eq:dom7} is of the form
\begin{equation}  \label{eq:dom71}  \begin{aligned} &
  - 2 \im e^{- \im xz_1}  m _{11}   \frac{\overline{ {c}}_1
   \overline{ m_{22}   }e^{ -\im x\overline{z}_1 -4\im t  \overline{{z}}_1^2} }{b^2} \times \\
   & \phantom{texttext} \left  (1
  +  O\left (  b^{-1}  {\left | c_1  m_{12}   e^{ \im xz_1 +4\im t {z}_1^2}     \right |}   \right )
  +     O\left (  b^{-1}  {\left |  m_{21} e^{  -\im x z_1  }     \right |}   \right )    \right )
 .
\end{aligned}\end{equation}
We claim  that the  quantity in \eqref{eq:dom71}  equals
\begin{equation}  \label{eq:dom9}  \begin{aligned} &
 - 2 \im \frac{e^{- \im xz_1}  \delta (z_1) ( \overline{\delta}  (z_1) )^{-1}
 \overline{ {c}}_1  e^{ -\im x\overline{z}_1 -4\im t  \overline{{z}}_1^2} }
{|e^{- \im xz_1}  \delta (z_1)    |^2+\left |\frac{  {{c}}_1   e^{  \im x {z}_1 +4\im t   {{z}}_1^2} }{  2 \im \Im  (z_1)}  {\delta  (z_1)}^{-1}  \right |^2 } (1+O(\epsilon t^{-\frac{1}{2}})) .
\end{aligned}\end{equation}
To prove this claim, we observe that since
$m_{ii}=\delta ^{-(-1)^{i^{}}} (z_1) +O(\epsilon t^{-\frac{1}{2}}) $  and
$|\delta ^{\pm 1} (z_1)|\ge \langle \rho \rangle ^{-2} $,
we have
\begin{equation*}   \begin{aligned} &
b^{ 2 }={|e^{- \im xz_1}  \delta (z_1)    |^2+
\left | \frac{c_1   e^{ \im x z_1 + 4\im t  z_1^2} }{  2 \im \Im  (z_1)}
{\delta  (z_1)}^{-1}  \right |^2 } (1+O(\epsilon t^{-\frac{1}{2}})).
\end{aligned}\end{equation*}
We have
$O\left (  b^{-1}  {\left |   c_1  m_{12}   e^{ \im xz_1 +4\im t {z}_1^2}   \right |}   \right )
= O(\epsilon t^{-\frac{1}{2}}) $    by
\begin{equation*}   \begin{aligned} &
      b^{-1}  {\left |  c_1  m_{12}   e^{ \im xz_1 +4\im t {z}_1^2}   \right |}  \le  \frac{\left |   m_{12}   e^{ \im xz_1 +4\im t  {z}_1^2}     \right |}{\left |   m_{22}   e^{ \im x z_1 + 4\im t z_1^2}     \right |}
 =  \frac{\left |   m_{12}       \right |}{\left |     m_{22}      \right |} \le C  \epsilon t^{-\frac{1}{2}}.
\end{aligned}\end{equation*}
We have  $  O\left (  b^{-1}  {\left |  m_{21} e^{  \im x\overline{z}_1  }     \right |}   \right )
= O(\epsilon t^{-\frac{1}{2}}) $    by
\begin{equation*}   \begin{aligned} &
      b^{-1}  {\left |    m_{21} e^{  -\im x z_1  }    \right |}
      \le  \frac{\left |   m_{21} e^{  -\im x z_1  }     \right |}{\left | m _{11}  e^{- \im xz_1}     \right |}
 =  \frac{\left |   m_{21}       \right |}{\left |     m_{11}      \right |} \le C  \epsilon t^{-\frac{1}{2}}.
\end{aligned}\end{equation*}
Hence \eqref{eq:dom9} is proved.

   Now we look at the term   in  \eqref{eq:dom9}.	
   For $z_1= \alpha _1+\im \beta _1$, $d_1= \log (\frac{|c_1|}{2 \beta _1} )$
   and  $\vartheta _1=\arg \left ( {  {c}_1 } \right )$,   dropping the
   factor $(1+O(\epsilon t^{-\frac{1}{2}}))$, for $\Delta   (z_1)$ defined in \eqref{eq:delta11}
   and inserting trivial factors $\Delta  /\Delta  =1$
	and  $\overline{\Delta}  /\overline{\Delta} =1$,
the expression in \eqref{eq:dom9}  equals
\begin{equation}  \label{eq:dom10}  \begin{aligned} &
   \frac{   -   4  \im    \beta _1 e^{-2\im \alpha _1x-4\im t (\alpha^2_1- \beta ^2_1)-\im \vartheta _1} \frac{\delta (z_1)  }{\Delta   (z_1)}     \frac{\overline{\Delta}  (z_1)}{\overline{\delta} (z_1)} \frac{\Delta   (z_1 )}{   \overline{\Delta}  (z_1)}   }
{e^{ 2\beta _1 x+8t\alpha _1\beta _1- d_1} |\frac{\delta (z_1)  }{\Delta   (z_1)}| \ |\Delta   (z_1)| +    e^{ -(2\beta _1 x+8t\alpha _1\beta _1- d_1)} |\frac{\Delta (z_1)  }{\delta   (z_1)}| \ |\Delta   (z_1)| ^{-1}}    .
\end{aligned}\end{equation}
Fix now a constant $\kappa >0$.  Then \eqref{eq:dom10} differs from the soliton solution
 \begin{equation}  \label{eq:dom11}
     -2\im \beta _1  e^{-2\im \alpha _1x-4\im t (\alpha^2_1- \beta ^2_1)-\im \vartheta _1 + 2\im \arg(\Delta   (z_1))}
 \text{sech}(2\beta _1 x+8t\alpha _1\beta _1- d_1+\log (|\Delta   (z_1)|))
\end{equation}
   by less than $c \kappa t^{-\frac{1}{2}} \epsilon $.  To prove this claim
   we observe that  the difference of \eqref{eq:dom10} and \eqref{eq:dom11}
   can be bounded, up to  a constant factor  $C=C(\omega _0, v_0)$,
		by the sum of the following two error terms:
   \begin{equation}  \label{eq:dom101}  \begin{aligned} &
   \frac{        \left |  \frac{\delta (z_1)  }{\Delta   (z_1)}     \frac{\overline{\Delta}  (z_1)}{\overline{\delta} (z_1)}  -1\right |   }
{e^{ 8 \beta _1 t|z_0- \alpha _1|   } (1+\| \widetilde{r} \| _{L^\infty(\R) } ^{2} )^{-1}}
\end{aligned}\end{equation}
and
\begin{equation}  \label{eq:dom102}  \begin{aligned} &
    \big |  \text{sech} \left (8 \beta t(-z_0+ \alpha _1) - d_1+\log (|\Delta   (z_1)|)\right )  \\& - \text{sech} \left (8 \beta _1t (-z_0+ \alpha _1) - d_1+\log (|\Delta   (z_1)|) +\log \left (\frac{|\delta   (z_1)|}{|\Delta   (z_1)|}\right ) \right  )   \big |
 .
\end{aligned}\end{equation}
We bound first \eqref{eq:dom101}.
For  $|z_0-\alpha _1|\ge  \kappa  t^{-\frac{1}{2}}$
formula \eqref{eq:dom101}  is bounded  by $ C   e^{- 8 \beta _1   \kappa\sqrt{t} }\epsilon$  by \eqref{eq:delta11}.  For  $|z_0-\alpha _1|\le  \kappa  t^{-\frac{1}{2}}$, for a fixed $K$ and
using  \eqref{eq:delta1}
 we bound    \eqref{eq:dom101} by
\begin{equation*}    \begin{aligned} &
      (1+\| \widetilde{r} \| _{L^\infty(\R) } ^{2} )  \left |  \frac{\delta (z_1)  }{\Delta   (z_1)}
			\frac{\overline{\Delta}  (z_1)}{\overline{\delta} (z_1)}  -1\right |   \le 4 \frac{C}{\sqrt{t}}
				\log (1+\| \widetilde{r} \| _{L^\infty(\R) } ^{2} )
     \le K   t^{-\frac{1}{2}}    \epsilon ^2.
\end{aligned}\end{equation*}
By Lagrange Theorem, \eqref{eq:dom102}  is bounded by
\begin{equation*}    \begin{aligned} &
   C
         \text{sech} \left (8 \beta _1t (-z_0+ \alpha _1) - d_1+\log (|\Delta   (z_1)|) +c\log \left (\frac{|\delta   (z_1)|}{|\Delta   (z_1)|}\right )  \right )     \left | \log \left (\frac{|\delta   (z_1)|}{|\Delta   (z_1)|}  \right )   \right |
\end{aligned}\end{equation*}
for some $c\in (0,1)$. This satisfies bounds similar to those
satisfied by  \eqref{eq:dom101}.

To complete  the proof of Theorem \ref{th:asstab} when   $u_0\in   H^1(\R) \cap L^{2,s}(\R)$,
we need to show that when one of  \eqref{eq:dom1}--\eqref{eq:dom5}  is false,
then the function in \eqref{eq:dom11}  is   $O(\epsilon t^{-\frac 12})$.
By Lemma \ref{lem:est m} the fact that  \eqref{eq:dom1}, resp.\eqref{eq:dom5},
false means  that for a fixed   $C=C(\rho _0)>0$ we have
\begin{equation*}
|   e^{-2\im  x z_1- 4\im  t \overline{z}_1^2}| =
   e^{   2(   \beta _1 x   +4 t\alpha _1\beta _1)}  \le C  \epsilon t^{-\frac 12}
\end{equation*}
and
\begin{equation*}
|   e^{  2\im  x \overline{z}_1+ 4\im  t  {z}_1^2}|  =
   e^{ -  2(   \beta _1 x   +4 t\alpha _1\beta _1)}  \le C  \epsilon t^{-\frac 12}.
\end{equation*}
Any of these  yields our claim that the function in \eqref{eq:dom11}  is   $O(\epsilon t^{-\frac 12})$.

This completes the proof of  Theorem \ref{th:asstab} for
$u_0\in   H^1(\R) \cap L^{2,s}(\R)$. Notice that
for $t\ge T( \varepsilon _0)$ the soliton in formula \eqref{eq:asstab2}
is given by formula  \eqref{eq:dom11}.

\vspace{0.2cm}

When  $   u_0\in   L^{2,s}(\R)  $ but $ u_0\not \in H^1(\R)$,
we consider a sequence $u_n \in H^1(\R) \cap L^{2,s}(\R)$ with $u_n\to u_0$ as $n \to \infty$ in $L^{2,s}(\R)$.
Then the sequence of spectral data from  $\{ u_n\} $ converges to the
spectral datum of $u_0$. This implies that for  $t\ge T( \varepsilon _0)$  we have \begin{equation}
\label{eq:asstab22}
\|  u_n(t,\cdot ) - \varphi _{\omega _n, \gamma _ +^{(n)} , v_n}(t,\cdot -x_+^{(n)} )\| _{L^\infty  (\R )}<
C \epsilon t^{-\frac{1}{2}},
\end{equation}
with a fixed constant $C$, since $C$ can be made to depend only on values of
$\varepsilon _0$ and $(\omega _0, v_0)$ in Theorem \ref{th:asstab}.
The sequence  $\{ (\omega _n  , v_n)\}$ converges to the parameters of the soliton
with spectral datum $(z_1,c_1)$ obtained  from the spectral datum     $(z_1,c_1,r)$ of  $u_0$. Finally,
$\{ (\gamma _ +^{(n)}  , x_+^{(n)})\}$ is a convergent sequence,
as can be seen in \eqref{eq:dom11} from their continuous dependence on the spectral data.
This means that for almost any $x$ and for any $t\ge T(   \varepsilon  _0)$, we have
\begin{equation*}
 \lim _{n\to \infty}
   \left ( u_n(t,x ) - \varphi _{\omega _n, \gamma _ +^{(n)} , v_n}(t,x -x_+^{(n)} ) \right )
   =  u(t,x ) - \varphi _{\omega _1, \gamma _ +  , v_1}(t,x -x_+ ).
\end{equation*}
Hence, bound \eqref{eq:asstab22} implies that for any $t\ge T(   \varepsilon  _0 )$, we have
\begin{equation*}
\|  u (t,\cdot ) - \varphi _{\omega _1, \gamma _ +  , v_1}(t,\cdot  -x_+ )\| _{L^\infty  (\R )}\le
C \epsilon t^{-\frac{1}{2}}.
\end{equation*}
The proof of Theorem \ref{th:asstab} is complete.

\qed

We end the paper  explaining the remark  that  the ground states  $\varphi_{\omega_1 ,\gamma_\pm ,v_1}(t,x-x_\pm )$ in the statement of Theorem \ref{th:asstab} are in general distinct.  The $+$  ground state   has been  computed explicitly in  \eqref{eq:dom11}.

\begin{lemma}
  \label{lem:distgs}   The  $-$  ground state   is given by formula  \eqref{eq:dom11} but with $\Delta   (z_1)$  replaced by
  \begin{equation} \label{eq:distgs1} \begin{aligned} &
     {\Lambda   (z_1) }  = \exp \left ( { { \frac 1{2\pi \im } \int _{\alpha _{1}}  ^{\infty}   \frac{\log (1+|{r} (\varsigma )|^2)}{\varsigma -z_1} d\varsigma} }\right ) .
\end{aligned}
\end{equation}
  \proof We  know that if $u(t,x)$ solves \eqref{eq:nls}  then $v(t,x) :=\overline{u}(-t,x) $
  solves the NLS with initial value  $ \overline{u}_0( x) $.  By standard arguments
  which can be derived from \eqref{eq:jost},  if $(r(z),  z_1 ,
 c_1 )$
   are the spectral data of   $  {u}_0 \in \mathcal{G}_1 $, then  we have
   $ \overline{u}_0 \in \mathcal{G}_1 $ with  spectral data  $(\overline{r}(-z),  -\overline{z}_1 ,
 -\overline{c}_1 )$. Using the latter,   by  \eqref{eq:dom11}   we then get for $ -t \nearrow \infty$
 \begin{equation*}  v(-t,x )\sim
      2\im \beta _1  e^{ 2\im \alpha _1x+4\im t (\alpha^2_1- \beta ^2_1)+\im \vartheta _1 - 2\im \arg( {\Lambda   (z_1)})}
 \text{sech}(2\beta _1 x+8t \alpha _1\beta _1- d_1+\log (| {\Lambda    (z_1)}|))
\end{equation*}
  with   $ \Lambda    (z_1)$ defined in terms of its complex conjugate
  (the following  is simply \eqref{eq:delta11} for the spectral data of   $ \overline{u}_0  $)
  \begin{equation*} \overline{\Lambda    (z_1)} :=
			 \exp \left ( \frac 1{2\pi \im } \int _{-\infty}^{ - \alpha _1}    \frac{\log (1+|{r} (-\varsigma )|^2)}{\varsigma +\alpha _1 -\im \beta _1} d\varsigma  \right )  .
\end{equation*}
Then  \eqref{eq:distgs1} is true.
Using $u(t,x)= \overline{v}(-t,x) $ and so taking the complex conjugate of the above formula, we obtain for
$  t \to - \infty$
 \begin{equation*}  u(t,x) \sim
      -2\im \beta _1  e^{ -2\im \alpha _1x-4\im t (\alpha^2_1- \beta ^2_1)-\im \vartheta _1 + 2\im \arg( {\Lambda   (z_1)})}
 \text{sech}(2\beta _1 x+8t \alpha _1\beta _1- d_1+\log (| {\Lambda    (z_1)}|))
\end{equation*}
thus completing the proof of  {Lemma}
  \ref{lem:distgs}.
\qed

\end{lemma}

\end{document}